\RequirePackage{fix-cm}
\documentclass[smallextended]{svjour3}       
\smartqed  
 \usepackage[text={7in,8.5in},centering]{geometry} 
\usepackage{graphicx}
\usepackage{amsmath, amssymb, amsfonts, amsgen,amstext,amsbsy, epsfig, graphicx,rotating, subfigure}
\usepackage[usenames]{color}
\usepackage{multirow,bm,xspace,booktabs,subfigure,enumerate,booktabs,xspace,setspace}
\usepackage{caption,tabulary,tabularx,ulem}
\RequirePackage[colorlinks,linkcolor=red,citecolor=blue,urlcolor=magenta]{hyperref}
\usepackage{paralist,mathtools,chngcntr}
\usepackage[titletoc,title]{appendix}
\usepackage[dvipsnames]{xcolor}
  \usepackage{epstopdf}
\usepackage{lineno}
\usepackage{natbib}
\usepackage{algorithm2e}
\RestyleAlgo{ruled}
\newcommand{\ie}{\mbox{\sl i.e.}}
\newcommand{\eg}{\mbox{\sl e.g.}}

\newcommand{\beq}{\begin{eqnarray*}}
\newcommand{\eeq}{\end{eqnarray*}}
\newcommand{\beqn}{\begin{eqnarray}}
\newcommand{\eeqn}{\end{eqnarray}}
\usepackage{amssymb,amsfonts,amsmath,breqn, color,float, url}
\usepackage{caption}
\usepackage{pdfpages}
\usepackage{pgfplots}
\pgfplotsset{compat=1.16}
\newcommand{\Xn}{X^{(N)}}

\newcommand{\Zn}{Z^{(N)}}
\newcommand{\Yn}{Y^{(N)}}

\newcommand{\zz}{{\mathbb Z}}
\newcommand{\nn}{{\mathbb N}}
\newcommand{\mm}{{\mathbb M}}
\newcommand{\rr}{{\mathbb R}}

\newcommand{\TT}{{\mathbb T}}

\newcommand{\veps}{\varepsilon}

\newcommand{\feq}{\end{eqnarray*}}
\newcommand{\feqn}{\end{eqnarray}}
\newcommand{\bec}{\begin{claim}}
\newcommand{\fec}{\end{claim}}
\newcommand{\becn}{\begin{claim*}}
\newcommand{\fecn}{\end{claim*}}

\begin{document}

\title{Deterministic Approximation of a Stochastic Imitation Dynamics with Memory
}
%



\author{Ozgur Aydogmus        \and
        Yun Kang 
}


\institute{O. Aydogmus \at
               Department of Economics, Social Sciences University of Ankara
                Ulus-Ankara, Turkey\\
              \email{ozgur.aydogmus@asbu.edu.tr}           
           \and
           Y. Kang \at
              Sciences and Mathematics Faculty
                College of Integrative Sciences and Arts, Arizona State University\\
              \email{yun.kang@asu.edu}
Mesa, AZ 85212, USA
}

\date{Received: date / Accepted: date}

\maketitle
\begin{abstract}
We provide results of a deterministic approximation  for non-Markovian stochastic processes modeling finite populations of individuals who recurrently play symmetric finite games and imitate each other according to payoffs. We show that a system of delay differential equations can be obtained as the deterministic approximation of such a non-Markovian process. We also show that if the initial states of stochastic process and the corresponding deterministic model are close enough, then  the trajectory of stochastic process stays close to that of the deterministic model up to any given finite time horizon with a probability exponentially approaching one as the population size increases.
We use this result to obtain that the lower bound of the population size  on the absorption time of the non-Markovian process is exponentially increasing. Additionally, we obtain the replicator equations with distributed and discrete delay terms as examples and analyze how the memory of individuals can affect the evolution of cooperation in a two-player symmetric Snow-drift game.  We investigate the stability of the evolutionary stable state of the game when agents have the memory of past population states, and implications of these results are given for the stochastic model.
\end{abstract}

\newpage

\section{Introduction}

Evolutionary game theory has been applied in economics, social and biological sciences where phenomena are typically aggregate outcomes of recurring strategic interactions in large populations of agents (see, {\it e.g.,} \cite{friedman1998economic,gintis2003explaining,broom2018biology,harms2008evolution}). In nature, those interactions are frequency-dependent that is the success of a player with a particular strategy depends on the number/frequency of agents adopting each strategy  (see,{\it e.g.,} \cite{aydogmus2020does,wakano2007social}). Even though there are some exceptions (see, among many, \cite{taylor2004evolutionary,hwang2013deterministic,aydogmus2017preservation,wang2017imitation}), the largest portion of models in the literature of evolutionary games hypothesize a continuum of interacting agents and identify the evolutionary processes as systems of ordinary differential equations  \citep{sandholm10}.  Each equation in a system keeps track of the population whose individuals are adopting one of the pure strategies. Replicator equations, for instance, can be considered as one of the most noteworthy examples of such models. For a wide class of game dynamics (including replicator equations), it was shown that evolutionary and non-cooperative games are strongly connected. In particular, the relationship between the steady states of the system of ODEs and the Nash equilibria of the game is given by the folk theorem of evolutionary game theory \citep[Theorem 1]{cressman2014replicator} (see also \cite{weibull1997evolutionary} for similar results). As an example, one may consider the result stating that time averages of trajectories of replicator equations approximate the mixed strategy Nash equilibrium of the game provided that there exists only one such Nash equilibrium \citep{hofsigm,weibull1997evolutionary}. This implies that even if a solution trajectory of the replicator equation fluctuates in time, the vector of time averages of the population shares constitutes a Nash equilibrium.

Regarding the relevance of the above-mentioned results, \cite{benaim} tried to answer an important question that is whether these deterministic models are good approximations of more realistic stochastic population processes used to model finite but large populations rather than infinite populations. Together a game played by a finite population of agents and a revision protocol that is used by agents to revise their strategies define a stochastic in particular, Markovian game dynamics \citep{sandholm10}.  Additionally, as noted by \cite{sandholm10}, these stochastic and deterministic models can be derived from a single foundation, {\it i.e.,} the above-mentioned deterministic evolutionary game models can be obtained as {\it fluid limits} or {\it mean-field equations} of the Markovian processes. Such a result can be proved by showing the probability that the trajectories of two processes (deterministic and stochastic) stay together up to some finite time $T$ approaches 1 as the population size goes to infinity (see, {\it e.g.,} \cite{sandholm10}). Yet these results are not sufficient to decide whether deterministic models are good approximations for their stochastic counterparts or not. The result obtained by \cite{benaim}, on the other hand, shows that the probability of large deviations of two processes exponentially approaches to zero. This result was classified as the strongest deterministic approximation result by \cite{sandholm10}. These approximations were used to obtain results regarding the behavior of the Markovian evolutionary processes. All of the above-mentioned results regarding the {\it fluid limits} and large deviation bounds are also valid for general Markov processes (not necessarily related to evolutionary games) and can be found in \cite{kurtz1970solutions, ethier} and \cite{norris}. Except for the above-mentioned studies, {\it mean-field} approximations for evolutionary games has been studied by many authors (see, {\it e.g.,} \cite{binmore97,binmore,boylan,borgers,corradi})

Here we would like to note that the above results are only valid if the Markov property ({\it i.e.,} memoryless agents) is assumed. This property implies that the payoff (or fitness) of a player at a certain time can be calculated if the population share of each pure strategy is known at that time. This assumption does not only imply that an agent ignores the recent population states when she updates her strategy but also means that she knows the exact frequencies of each phenotype at the time of the update. The assumption of memoryless agents might be more relevant in a biological setting. For an evolutionary birth-death process in which the players' reproduction rate depends on the fitness of each phenotype as considered by  \cite{taylor2004evolutionary} and \cite{fudenberg2006evolutionary}, this assumption makes more sense, since agents do not need to calculate their fitness and decide accordingly in such a scenario. The prisoners' dilemma game played by RNA viruses can be considered as an example of this biological setting \citep{turner1999prisoner}. Considering individuals with cognitive abilities deciding how to update their strategies, on the other hand, requires to relax this assumption and consider time-delayed and/or averaged information acquisition into account. 

This problem has been pointed out by many authors (see, {\it e.g.,} \cite{alboszta2004stability,miekisz2011stochasticity,moreira2012evolutionary,wang2017imitation,yi1997effect,bodnar2020three}) studying deterministic evolutionary dynamics in infinite populations via delay differential/difference equations. In particular, \cite{mikekisz2008evolutionary} notes the importance of considering delayed information as follows:
\begin{quote}
``{\it It is very natural, and in fact important, to introduce a time delay in the
population dynamics; a time delay between acquiring information and acting
upon this knowledge or a time delay between playing games and receiving payoffs.}"\end{quote} The above explanation is more relevant when the imitation dynamics between individuals with cognitive abilities is considered. \cite{wang2017imitation} studied such deterministic imitation dynamics and suggested studying the stochastic imitation dynamic for finite populations.

In this study, we consider a general class of revision protocols that takes the history of the population states into account in modeling evolution in large but finite populations. Specifically, a non-Markovian process considered here is a direct generalization of the model studied by \cite{benaim} relaxing the assumption of the memoryless agents. So we consider a population of $N$ agents adopting one of the pure strategies in a d-player normal form game. The evolution takes place by allowing only one player to revise her strategy at times labeled by $0,\delta,2\delta,\cdots$ where $\delta=N^{-1}.$ Agents revise their strategies according to the payoffs they obtain, {\it i.e.,} by combining their knowledge of payoff matrix and the state of the population. The state of this population process is determined by the set of $m$ consecutive $d$-dimensional vectors in the $d-1$ dimensional unit simplex denoting the history of population shares for the last $m$ updates. The state of the process may be used to calculate the payoff or fitness of an individual by considering a discrete-time delay ({\it i.e.,} the strategy revision at time $t$ depends on the population shares at time $t-m\delta$)  or a distributed time delay ({\it i.e.,} the strategy revision at time $t$ depends on a weighted average of population shares of these $m$ consecutive vectors).

We extend the deterministic approximation results for the Markov processes via ODEs by showing that {\it fluid limits} of above-described population processes are delay differential equations and that the probability of large deviations of trajectories of two processes is also exponentially bounded above. To the best of our knowledge, this is the first study deriving such a large deviation bound for non-Markovian processes. Using these approximations we also obtained an exponentially increasing lower bound in population size for the absorption (or fixation) times when the trajectory of the delay differential equations is bounded away from the boundary of the unit simplex. We obtain replicator equations with discrete and distributed delays from a microscopic update rule taking the history of the process into account. We show delayed replicator equations also satisfy time averaging property, and an extension of this result for the stochastic model is also given.

The text is organized as follows: In Section \ref{modell}, the notation and the model will be introduced. In Section \ref{ffipba}, the deterministic approximation results and their implications for the stochastic process in terms of absorption times will be given. In Section \ref{dreww}, delayed replicator equations will be obtained and implications of our results will be given for the snow-drift game. We conclude the paper in Section \ref{conc} and  defer the proofs to Section \ref{prfs}.


\section{Notation and the Model}
\label{modell}
Let $d\geq 2$ be a fixed integer and introduce the vector notation $\mathbf x=(x_1,x_2,\cdots,x_d)\in\rr^d.$ We start by considering a symmetric two-person games with a pure strategy set $S=\{1,2,...,d\}$ and the mixed strategy simplex
\beq
\Delta_d=\Big\{\mathbf x\in\rr_+^d\,\Big|\,\sum_{i\in S} x_i=1\Big\},
\feq
where $\rr_+$ is used to denote the set of non-negative real numbers $\{y\in\rr\,|\,y\geq 0\}.$ The pure strategies of the game are identified with corners of the simplex. In particular, $i\in S$ is identified by the unit vector $\mathbf e_i=(0,\ldots,0,1,0,\ldots,0)\in \Delta_d$ whose only non-zero
component is $1$ at the $i^{\text{th}}$ place. 
\par
 Throughout the paper we use the maximum ($L^\infty$) norms for vectors in $\Delta_d$ and functions in $C(\Delta_d),$ the space of real-valued continuous functions on $\Delta_d$ that are defined as follows:
\beq
\|\mathbf x\|:=\max_{i\in S}|x_i|~\mbox{\rm for}~\mathbf x\in \Delta_d
\qquad \mbox{\rm and}\qquad 
\|f\|:=\max_{\mathbf x\in \Delta_d}|f(\mathbf x)|~\mbox{\rm for}~f\in C(\Delta_d).
\feq In addition to these norms, we frequently employ the norm defined on the function space  $C:=C\bigl([-r,0],\Delta_d\bigr)$ defined as follows: \beq \|\boldsymbol\phi\|_C=\sup\bigl\{\|\boldsymbol\phi(\theta)\|~:\,-r\leq\theta\leq0\bigr\} \mbox{ for any } \boldsymbol\phi\in C\eeq

In the following lines, we describe our stochastic and deterministic models, and give some of their properties.


\subsection{Stochastic model}\label{stm}
Stochastic processes we study in this paper describe the evolution in time of a population of $N$ individuals. Here $N\in\nn$ is a constant number, and we are mostly interested in the behavior of the system for large but finite values of $N.$  

Suppose that each individual is associated with a pure strategy in $S$ at any instance of time and the players' strategies evolve in time stochastically. Rules of this evolution are determined through an interaction between the players as specified in the following lines.  An individual in the population is said to be an $i^{th}$-strategist if she is presently associated with the pure strategy $i.$ Only at times $\tau\in \TT=\{0,\delta,2\delta,\cdots\},$ where $\delta=1/N,$ exactly one randomly chosen individual is given an opportunity to change her strategy.  
\par
Specifically, we will concentrate on the sequence of $d$-dimensional vectors $\mathbf\Xn(\tau)=\bigl(\Xn_1(\tau),\Xn_2(\tau),\ldots,\Xn_d(\tau)\bigr),$ $\tau\in\TT,$ with
\beq
\Xn_i(\tau)=\frac{\#\{i^{th}~\mbox{\rm strategists at time $\tau$}\}}{N}, \qquad i\in S.
\feq
We refer to $\mathbf\Xn(\tau)$ as the {\it population profile at time $\tau.$} To define the state space of the process, we need the following set:
\beq
\Delta_{d,N}=\Bigl\{0,\frac{1}{N},\frac{2}{N},...,\frac{N-1}{N},1\Bigr\}^d.
\feq
Using this we define the state space: $\mathbf\Xn(\tau)\in\Delta_d^N:=\Delta_d\cap\Delta_{d,N}$ for all $\tau\in\TT.$ The specific model to be studied here forms a not-necessarily Markov process yet it is a generalization of the Markov chain studied by \cite{benaim}. For any pair $i,j \in S$, the transition probabilities of the aforementioned Markov process is determined by a function $\tilde p_{ij}:\Delta_d\to[0,1]$ satisfying $\tilde p_{ij}(\mathbf x)=0$ if $x_j=0$ and defined as\beq\tilde p_{ij}=Pr\Big[ \mathbf\Xn(\tau+\delta)=\mathbf x+\delta(\mathbf e_i-\mathbf e_j)\Big|\mathbf\Xn(\tau)=\mathbf x \Big]\eeq \cite{benaim} supposed that the conditional probability that a $j^{th}$-strategist will become an $i^{th}$-strategist (and hence $\tilde p_{ij}$) is continuous in the current state $x.$

Employing the above-given memoryless process in a social model implies that the agents can collect the data regarding the frequencies (or numbers) of each type of agent $\mathbf\Xn(\tau)$  and process this information instantaneously to calculate the imitation probabilities. Here, we assume that these individuals imitate another agent's strategy with a probability depending on payoffs of each type calculated using past frequency vector $\mathbf\Xn(\tau-r)$ for some $r>0$ or some weighted average of these vectors. Here $r=m\delta$ for some constant $m\in \nn.$ Hence we define the following set of integers: $\mm =\{k\in\zz~:\,-m\leq k\leq 0 \}.$
\par

To incorporate a time delay in our model, we start by considering the $\sigma$-algebra $\mathcal F_\tau, ~~\tau\in\TT$ generated by \beq\big\{\mathbf\Xn(k\delta)=\boldsymbol\xi_k~\big|~k\in\mm\big\}\cup\big\{\mathbf\Xn(\kappa\delta)~\big|~0< \kappa\leq \tau N\big\}\feq for vectors  $\boldsymbol\xi_k\in \rr^d.$
We use this filtration to define the following set indicating the difference of two $\sigma$-algebras: \beq \mathcal A_\tau=\mathcal F_\tau\backslash\mathcal F_{\tau-r-\delta},~~~\tau\in\TT.\eeq $A_\tau$ contains all the information necessary to determine the solution to the stochastic process $\mathbf\Xn(\sigma)$ for all $\sigma>\tau ~(\sigma,\tau\in \TT).$

Now we consider the following conditional probability \beq p_{ij}^m=Pr\Big[ \mathbf\Xn(\tau+\delta)=\mathbf\Xn(\tau)+\delta(\mathbf e_i-\mathbf e_j)\Big|\mathcal A_\tau\Big]\eeq   satisfying $ p_{ij}^m=0$ if $\Xn_j(\tau)=0.$ This implies that the probability that an $i^{th}$-strategist imitates a $j^{th}$-strategist is positive only if there exists at least one $j^{th}$-strategist in the population.  Then the transition probabilities of the stochastic process for any $\mathbf v\in \rr^d$ are \beqn\label{tr} Pr\Big[ \mathbf\Xn(\tau+\delta)=\mathbf\Xn(\tau)+\delta \mathbf v\Big|\mathcal A_\tau\Big]=\begin{cases} p_{ij}^m, &\mbox{if } \mathbf v=\mathbf e_i-\mathbf e_j \\
0,&\mbox{otherwise }.\end{cases}\eeqn Here we would like to note that above defined transition probabilities reduce to that of the Markov process  studied by \cite{benaim} if $m=0.$ If, on the other hand, $m\in\zz_+,$ then the stochastic process is no more Markovian and it can be used to study the effects of information delays on the game dynamics. In the latter case, the conditional probabilities $p_{ij}$ are not simply a function taking values from $\Delta_d$ and having images in $[0,1],$ since the initial condition depends on vectors $\boldsymbol\xi_k \in\mathcal A_0.$

To obtain the mean-field equations of this process in the following section, we need to specify probabilities $p_{ij}^m$ given in \eqref{tr}. It is clear that $p_{ij}^m$ and hence the transition probability at time $\tau\in \TT$ depends on the values of $\Xn(\tau+k\delta)$ for all $-m\leq k\leq 0.$ Hence we denote this quantity as follows: \beqn\label{discretepr} p_{ij}^m=p_{ij}^m\Bigl(\mathbf{X}_{\boldsymbol\tau}^{\mathbf{(N)}} \Bigr)\eeqn where \beqn \mathbf{X}_{\boldsymbol\tau}^{\mathbf{(N)}}(k\delta)=\mathbf{\Xn}(\tau+k\delta),~ k\in\mathbb M=\{-m,-m+1,\cdots,0\}.\label{notds}\eeqn Clearly for $m=0,$ we have  $p_{ij}^0=\tilde p_{ij}.$ If the transition probabilities depend only on the discrete time delays then the formulation of $p_{ij}^m$ given in \eqref{discretepr} would be enough to obtain the mean-field equations. To get results concerning models taking more general delays into account, we need to consider a continuous time version of $p_{ij}^m$ defined for any $t\in\rr_+$ as follows: \beqn\label{continuouspr} p_{ij}\Bigl( \overline{\mathbf{X}}^{\mathbf{(N)}}_{\mathbf t}\Bigr):=p_{ij}^m\Bigl(\mathbf{X}_{\boldsymbol\tau}^{\mathbf{(N)}}\Big) \mbox{ for } t\in[\tau,\tau+\delta)\eeqn where $\overline{\mathbf{X}}^{\mathbf{(N)}}_{\mathbf t}$ is defined for any $\theta\in[-r,0]$ as follows: \beqn \overline{\mathbf{X}}^{\mathbf{(N)}}_{\mathbf t}(\theta)=\mathbf \Xn(\sigma)\mbox{ for } t+\theta\in[\sigma,\sigma+\delta).\label{notcs}\eeqn This notation is borrowed from the delay differential equations literature and will be used frequently in the remainder of the paper.

Finally we would like to state the following assumption regarding the vectors $\boldsymbol\xi_k:$ \begin{itemize}
    \item [$A$-] $\|\boldsymbol{\xi_{k}}-\boldsymbol{\xi_{k+1}}\|\,\leq 2\delta$ for any $k\in\mm.$
\end{itemize} Recall the assumption that exactly one randomly chosen individual is given an opportunity to chance her strategy at times $\tau\in \TT,$ which implies that  $\|\mathbf{\Xn}(\tau)-\mathbf{\Xn}(\tau+\delta)\|\,\leq2\delta.$ Hence Assupmtion $A$ is a natural extension of this property for the history of the process $\mathbf\Xn(k\delta)$ for $k\in\mm.$


\subsection{Mean-field equations and semiflows}\label{mnf}
Hereafter, we assume that $m\in \zz_+$ or $r>0.$  Our aim is, now, to obtain the mean field equations of the above given process. To do this we need to determine $ N\cdot E\Big[{\mathbf{\Xn_i}}(\tau+\delta)-\mathbf{\Xn_i}(\tau)\Big|\mathcal A_\tau\Big]$ ($\ie$ the expected net increase in the number of $i^{th}$-strategists from one transition time to next conditioned on $\mathcal A_\tau$) which is given by \beqn\label{veca} F_i (\mathbf x_t):=\sum_{k\neq i} p_{ik}(\mathbf x_t) -\sum_{k\neq i}p_{ki}(\mathbf x_t)\eeqn where $\mathbf x_t$ is defined by \beq \mathbf x_t(\theta):=\mathbf x(t+\theta),~~ -r\leq\theta\leq0.\eeq Hence the associated mean-field equations are \beqn\label{delayde} \dot x_i=F_i(\mathbf x_t),~~i\in S, \mathbf x_t\in C\\\mathbf x_0(\theta)=\phi(\theta),~~ -r\leq\theta\leq0 \nonumber \eeqn where $\dot x_i$ is used to denote the derivative of $x_i$ with respect to $t$ and $\phi\in C$ is the initial function.

We assume that $F:C\to R^d$ satisfies the following Lipschitz condition which is required for the existence and uniqueness of solutions: 
\begin{itemize}
    \item[($Lip$)] For any $M>0,$ there exists a $K>0$ such that \beq\|\mathbf F(\boldsymbol\phi)-\mathbf F(\boldsymbol\psi)\|\leq K\|\boldsymbol\phi-\boldsymbol\psi\|_C,\mbox{ for any } \|\boldsymbol\phi\|_C,\|\boldsymbol\psi\|_C\leq M\eeq  
\end{itemize} where $\|\boldsymbol\phi\|_C$ is the sup norm defined in $C.$

We have the following result on the existence, uniqueness, and continuation of solutions to \eqref{delayde} which follows from  \citet[Theorem 3.7 and Proposition 3.10]{smith2011} 

\begin{lemma}\label{exuniq} Suppose that $\mathbf F$ is continuous and satisfies the Lipschitz condition ($Lip$) for some $M>0.$ If $\|\boldsymbol\phi\|_C\leq M$ then there exists a unique solution $\mathbf x(t)=\mathbf x(t,\boldsymbol\phi)$ of \eqref{delayde} defined for all $t>0.$
\end{lemma} A sketch of the proof is given in Section \ref{pexuniq} and it basically follows \cite{smith2011}.

 When dealing with ordinary differential equations (or Markov processes), the state of the ODE corresponds to its solution. For delay differential equations, we no longer have this flexibility. The solution through $(t,\phi)$ is denoted by $\mathbf x(t,\phi).$  We, on the other hand, denote the state of delayed system \eqref{delayde} by $\mathbf x_t$ which contains all the necessary information to determine $\mathbf x(s)=\mathbf x(s,\boldsymbol\phi)$ for $s>t.$ In particular, we have \beq\mathbf x_t(\boldsymbol\phi)(\theta)= \mathbf x(t+\theta,\boldsymbol\phi)\eeq for $\theta\in [-r,0].$ Using the fact that the system is autonomous, the system of equations \eqref{delayde} defines a semiflow as follows:
\beq \boldsymbol\Phi(t,\boldsymbol\phi)=x_t(\boldsymbol\phi)=S(t) \boldsymbol\phi\eeq satisfying $S(0)\boldsymbol\phi=\boldsymbol\phi$ and $S(t)S(s)\boldsymbol\phi=S(t+s)\boldsymbol\phi$ for $t,s\geq0$ and $\boldsymbol\phi\in C$ (see, $\eg,$ \cite{smith2011,sell2013dynamics}).   Hence the positive trajectory through $\boldsymbol\phi$ is defined as \beq \gamma^+(\boldsymbol\phi)=\{S(t)\boldsymbol\phi~:\,t\geq0\}.\eeq Similarly, for any set $B\subset C,$ we define trajectories through $B$ as \beq \gamma^+(B)=\{S(t)u~:\,u\in B,\,t\geq0\}.\eeq Here we say the set $A$ is invariant   if $\gamma^+(A)=A.$

 We say that the set $A$ attracts $B$ if for any $\veps>0$ there exists a $\mathcal T=\mathcal T(\veps)$ such that \beq d_C(S(t)B,A)\leq \veps \mbox{ for any } t\geq\mathcal T\eeq where  $d_C$ is the metric induced by the supremum norm of $C.$

Hence, we say $\mathcal A$ is an attractor provided that \begin{itemize}
    \item $\mathcal A$ is a compact invariant set in $C;$
    \item There is a neighborhood $U\subset C$ of $\mathcal{A}$ such that $\mathcal A$ attracts every bounded set in $U.$
\end{itemize} Lastly, we define the basin of attraction of an attractor $\mathcal A$ as follows: \beq B(\mathcal A)=\bigl\{ \boldsymbol\phi \in C~:\, d_C\bigl(S(t)\boldsymbol\phi,\mathcal A\big)\to0 \mbox{ as } t\to\infty \bigr\}.\eeq The above given definitions regarding semiflows can  be found in \cite{kuang} and \cite{sell2013dynamics}. 


\section{From Finite to Infinite Populations and Back Again}
\label{ffipba}
In this section, we aim to understand how the non-Markovian stochastic process and its mean-field equations are linked. In the literature, such approximation results are well-established for processes having Markovian property and ordinary differential equations \citep{ethier,kurtz,norris}. In particular, these results were used to reveal the relationship between Markovian game dynamics in finite populations and their mean-field equations (see, $\eg,$ \cite{benaim,binmore97,binmore,boylan,borgers,corradi,sandholm10}). To the best of our knowledge, deterministic approximations to non-Markovian stochastic processes have not been studied. 

\subsection{Delay differential equation approximation}
Here we aim to find a heuristic law of large numbers result stating that the trajectories of the non-Markovian stochastic process with high probability stay in close proximity of solutions to the associated deterministic delay equations during any given bounded time interval, provided that the population is large enough. A similar heuristic law of large numbers result showing the link between Markovian population game models and their mean-field equations is given by \cite{benaim}, and \cite{sandholm10} states that this result is the strongest deterministic approximation result in the literature.  

To measure the fit of the deterministic approximations over bounded time intervals, we rewrite the the stochastic process \eqref{tr} in continuous time by defining the {\it interpolated process} as follows: 
\beqn
\label{interpola}
\mathbf\Yn(t)=\mathbf\Xn(\tau)+\frac{t-\tau}{\delta}\bigl(\mathbf\Xn (\tau+\delta)-\mathbf\Xn (\tau)\bigr) \qquad \forall\,t\in \bigl[\tau,\tau+\delta\bigr)
\feqn which is defined for all $\tau\in \{-m\delta, (1-m)\delta,\cdots,-\delta\}\cup\TT.$

Using this process, we define the stochastic variable describing the maximal deviation in any population share on the bounded time interval $I:=[0,T]$ as follows: \beqn\label{devia} D\bigl(T,\boldsymbol{\phi}\bigr)=\max_{t\in I}\bigl\|\mathbf\Yn(t)-\mathbf x(t,\boldsymbol{\phi}))\bigr\|.\eeqn Regarding this stochastic variable we have the following result:

\begin{theorem}\label{mainth}
Suppose that Assumption $A$ holds and the initial function $\boldsymbol{\phi}\in C$ satisfies  $\phi(k\delta)=\boldsymbol
\xi_k$ for all $k\in \mm.$ \\
Then, for any $\veps>0,$ $T>0$ and large enough population size $N\in \nn,$ there exists a constant $c>0$ such that \beq Pr\Bigl[D\bigl(T,\boldsymbol{\phi}\bigr)\geq\veps~\big|\,\mathbf{\Xn_0}(k\delta)=\boldsymbol{\phi}(k\delta),\, k\in\mm\Bigr]\leq 2d \exp(-\veps^2 c N)\eeq
\end{theorem} The proof of this theorem, which is an extension of the result obtained by \cite{benaim} for Markov processes, is in Section \ref{pmainth}.

\subsection{Results on Absorption Times}
\label{exits}
Here, we use Theorem \ref{mainth} to obtain results regarding absorption times. Such results heavily depend on the results regarding the exit times from subsets of $\Delta_d$ (see, Section \ref{pth3}). 
  
  Yet we need to make sure our process modeling the imitation dynamics has absorbing states. Following the discussion on metastability in \cite[pp. 885-886]{benaim}, we assume the followings to guarantee existence of absorbing states:
  \begin{itemize}
      \item[C1-] $X^{(N)}_k(\tau)\in(0,1)$ implies $p^m_{kj}\Bigl(\mathbf{X}_{\boldsymbol\tau}^{\mathbf{(N)}}\Bigr)>0$ for some $k\neq j.$
      \item[C2-] $X^{(N)}_j(\tau)=0$ implies $p^m_{kj}\Bigl(\mathbf{X}_{\boldsymbol\tau}^{\mathbf{(N)}}\Bigr)=0$ for all $k\neq j.$
  \end{itemize} C1 implies that if some but not all individuals in the population uses $k^{th}$ strategy at time $\tau$ then a $k^{th}$ strategist is able to adopt another strategy at time $\tau+\delta.$ C2, on the other hand, implies that probability that a $k^{th}$ strategist change her strategy to $j$ is zero at time $\tau+\delta$ if there is no $j^{th}$ strategist in the population at time $\tau.$ As pointed out in \cite[Remark 2]{benaim}, if both of these conditions are satisfied the stochastic process reaches the boundary of the simplex and stays there forever. Hence, the population distribution stays put if all individuals adopt the same strategy. This phenomenon is called as the fixation of the population in the literature. An important quantity of interest in the dynamics of finite populations is the average time until fixation occurs (see e.g., \cite{traulsen2009stochastic,ewens2004mathematical}).
  
  Through this section, we have the following assumption regarding the attractor $\mathcal A:$
  \begin{itemize}
      \item[A-]  $\mathcal A\subset \tilde C:=C\bigl([-r,0], D\bigr)$ where $D\subset \mbox{int}(\Delta_d)$ is a closed set. Here, $\mbox{int}(\cdot)$ is used to denote the interior of a set. This implies that $\boldsymbol\phi(\theta)\neq \mathbf e_i$ for any $i\in S$ and $\theta \in [-r,0]$ provided that $\boldsymbol\phi\in \mathcal A.$
  \end{itemize}

  This assumption states that functions contained in the attractor of the semi-dynamical system never touches the boundaries of $\Delta_d.$ This guarantees that the trajectory of the deterministic process stays away from the boundary of the unit simplex given the initial condition. In this case, Theorem \ref{mainth} states that the sample trajectories of the stochastic process should be close to that of the deterministic process which implies that the sample paths should be away from the boundaries of the simplex containing the absorbing states. To study this phenomenon analytically, define the absorption/fixation time as follows:   \beq T_a^N=\inf\{\tau\in \TT~:\,\mathbf{X}^{\mathbf{(N)}}(\tau)=\mathbf e_i \mbox{ for some } i\in S \} \eeq Now we state our first result concerning the mean absorption probability.

\begin{corollary}\label{corr} Suppose that $\mathbf Y^N_0\in B(\mathcal A)$ for all $N.$ Then there exists a constant $\alpha>0$ such that

  \beq E\bigl[T_a^N\bigr]\geq\frac{1}{4d}e^{\alpha N}-1.\eeq
\end{corollary}

The proof of this result is given in Section \ref{pth31}. This result tells us that the mean time to absorption exponentially increases with the population size $N.$ Similar results are given for Markovian evolutionary games by \cite{benaim} and for extinction time of an epidemic model by \cite{aydogmus2016extinction}.  In addition, results obtained in Theorem \ref{mainth} allow us to utilize Borel-Cantelli Lemma to show that the absorption time exceeds any upper bound as the population size goes to infinity.

\begin{corollary}\label{cor2}
\label{liminf}
Suppose that $\mathbf Y^N_0\in B(\mathcal A)$ for all $N.$ Then we have
\beq
Pr\Bigl[\liminf_{N\to\infty}T^{N}_a=+\infty\Bigr]=1.
\feq
\end{corollary} A sketch of the proof of this assertion is given in Section \ref{pcor2}.


\section{Delayed Replicator Equations}
\label{dreww}
\subsection{Replicator Equations with Distributed and Discrete Delays}
Here our aim is to obtain delayed mean-field equations for a known imitation rule, {\it i.e.,} replicator rule. We, first, determine $p_{ij}^m$ which is the probability that an $i^{\mbox{th}}$ strategist becomes $j^{\mbox{th}}$ strategist. Suppose that, at time $t,$ an $i^{\mbox{th}}$ strategist is chosen with probability $x_i(t)$ who choose a $j^{\mbox{th}}$ strategist to imitate with probability $x_j(t).$ As in the case of Markovian setting, we suppose that the former individual imitates the latter with a probability proportional to $\big[f_i(\mathbf x_t)-f_j(\mathbf x_t)\big]_+$ where $f_i(\mathbf x_t)$ is the history dependent fitness of the $i^{\mbox{th}}$ strategists. Using this rule an agent imitates the opponent only if
the opponent’s payoff is higher than her own. Thus we have \beq p_{ij}~ \sim \,x_i(t)x_j(t)\big[f_i(\mathbf x_t)-f_j(\mathbf x_t)\big]_+.\eeq Hence, the delayed vector field \eqref{veca} for this type of comparison rules can be determined as follows \beqn\label{dreplicator} F_i (\mathbf x_t):=x_i(t)\Big(f_i(\mathbf x_t)-\sum_{k\in S} x_k(t)f_k(\mathbf x_t)\Big).\eeqn With this nonlinear function, \eqref{delayde} determines the delayed replicator equation. The simplest form of this equation can be determined by considering a single discrete delay as follows: \beq f_i(\boldsymbol{\phi})=\mathbf e_i'A\boldsymbol{\phi}(-r)\eeq where $\boldsymbol\phi\in C,$ $A\in \rr^{d\times d}$ is the payoff matrix and $\mathbf e_i'$ is the transpose of the vector $\mathbf e_i.$ With this function, the replicator equation with a single discrete delay is given by \beqn\label{ddreplicator} \mathbf{\dot x}(t)=\mathbf x(t)\cdot\Big( A\mathbf x(t-r)- \mathbf x'(t) A\mathbf x(t-r)\Big)\eeqn where $\cdot$ is used to denote the component wise (or Hadamart) product of two vectors and $\mathbf x'$ is used to denote the transpose of the vector $\mathbf x.$ 

We would like to note that it is also possible to get a replicator equation using distributed time delays. In particular consider a probability measure $\mu$ satisfying $\int_{-r}^0\,d\mu(s)=1,$ and let $\mathbf{\Bar{x}}(t)=\int_{-r}^0\mathbf{x}(t+s)\,d\mu(s).$ Then the replicator equation with distributed time delays is given by \beqn\label{distdreplicator} \mathbf{\dot x}(t)=\mathbf x(t)\cdot\Big( A\mathbf{\Bar{x}}- \mathbf x'(t) A\mathbf{\Bar{x}}\Big).\eeqn We would like to note that both of these delayed replicator equations reduce to replicator ODE for $r=0.$ For both of these models \eqref{ddreplicator} and \eqref{distdreplicator}, we have the following time averaging property: 

\begin{theorem}\label{timeav} Let $x(t)$ be a solution to one of the delayed replicator equations \eqref{ddreplicator} or \eqref{distdreplicator}. 
\begin{itemize}
    \item[a] If there is no interior equilibrium of replicator ODE, then $x(t)$ approaches the boundary of the simplex $\Delta_d$ asymptotically.
    \item[b] If there exists a unique interior equilibrium $\mathbf p$ of the replicator ODE, then \beq p_i=\lim_{T\to\infty}\frac{1}{T}\int_0^Tx_i(t)\,dt\feq  provided that the solution $\mathbf x(t)$ in the interior of the simplex $\Delta_d$ remains bounded away from the boundary of the simplex.
\end{itemize}
\end{theorem} The proof of this result is in Section \ref{ptimeav}. The above-given result is well-known for the replicator ODEs. We show that it is also valid for delayed replicator equations.  A corollary of this is given as follows:

\begin{corollary}\label{cor121} Suppose that $\mathbf p$ is the unique interior rest point of the replicator ODE.  Let $B_\veps(\mathbf p)$ be a ball with radius $\veps>0$ centered at $\mathbf p$ then for any $\veps>0$ there exists $\tau_\veps\in\TT$ such that
\beq
Pr\Biggl( \frac{1}{\tau N}\sum_{0\leq\sigma\leq\tau} \mathbf{X}^{\mathbf{(N)}}({\sigma})\not\in B_\veps(\mathbf p)\Bigg)\leq 2d\exp\Bigl(-\veps^2 c_\tau N\Bigr)
\feq for any $\tau>\tau_\veps.$
\end{corollary}
In Section \ref{pc121} we give a proof of this corollary which is a simple modification of the proof of the Theorem \ref{mainth}.

\subsection{Hawk-Dove Game as an Example}
\label{sec:11}
Here we consider a two-player game with strategies $A$ and $B$ with payoff matrix: 
\[\begin{tabular}{ l | c  r }
     & A & B \\ \hline
    A & $a$ & $ b$ \\ 
    B & $c$ & $ d$ \\
  \end{tabular}
\]

To describe the fitness of each phenotype, consider a population of $N$ individuals. If the frequency of type A individuals at time $\tau\in \TT$ is denoted by $\Zn(\tau),$ then the average frequency of this type is given by \beq\bar Z^{(N)}(\tau)=\sum_{i=0}^m k_i \Zn(\tau-i\delta)\eeq where $\mathbf k=(k_0,k_1,\cdots,k_m)$ is a discrete probability distribution.  Using this average, we calculate the payoffs of these two phenotypes as follows: \beq f_A=a\bar Z^{(N)}(\tau)+b\bigl(1-\bar Z^{(N)}(\tau)\bigr) \text{  and  }  f_B=c \bar Z^{(N)}(\tau)+d\bigl(1- \bar Z^{(N)}(\tau)\bigr).\eeq

If two individuals are randomly chosen from the population  at each time in $\TT,$ then the following actions take place:  One of these individuals is chosen as a role model and the other (focal) individual adopts the strategy of the role model with a probability depending on the payoff difference.

In particular, the specific form of the probability that a randomly chosen individual with strategy B adopts strategy A is given by  
\beq p_{AB}^m=\Zn(\tau)\bigl(1-\Zn(\tau)\bigr)\bigl[f_A-f_B\bigr]_+.\eeq  We would like to note that this probability is the imitation rule taking time delays into account (see, $\eg$ \cite{hofsigm}).

The above-mentioned process describes a Markov-chain for $m=0$. However, in
reality, it is not easy to imitate the role model's
strategy according to the comparison between the current
payoffs to her own and to role model's strategies, $\ie,$ calculating $f_A$ and $f_B$ for $m=0$ requires the knowledge of the immediate frequency of each type. In particular, calculating such payoffs requires collecting the data regarding the frequencies (or numbers) of each type of agents $\Zn(t)$ instantaneously. A more reasonable assumption is that an individual imitates its role model's strategy with a probability depending on past payoffs of focal and role model agents' strategies. Such an assumption leads us to investigate the effect of time delay in imitation dynamics.

Using the probability function $p_{AB}^m$ we obtain the following delayed replicator equations as the mean-field equations:
\beqn\label{delrep} \dot z=z\bigl(1- z\bigr) \Bigl((a-b-c+d)\bar z+b-d \Bigr)\eeqn where $z:=z(t)$ is the frequency of type A individuals at time $t$ and $\bar z:=\int_{-r}^0z(t+s)K(s)\,ds$ is the average frequency of individuals of type A. Here, note that the discrete probability distribution $\mathbf k$ is a discretization of the probability kernel function $K(s)$ satisfying $\int_{-r}^0K(s)\,ds=1.$ 

Now we would like to study how this delay term effects the dynamics. In particular, we linearize the delayed replicator equation by taking $z(t)=e+\veps\zeta(t) $ for small $\veps>0,$ where $e$ is an equilibrium point of \eqref{delrep}. Plugging this ansatz into the equation, at level $O(\veps),$ we obtain \beq \dot \zeta=(b-d)\zeta ~\mbox{ and } ~\dot \zeta=(a-c)\zeta \eeq for $e=0,1,$ respectively. Hence the delayed term does not affect the stability of the pure strategy equilibria 0 and 1. On the other hand, if $e\in(0,1)$ is an interior equilibrium, then we have \beq \dot \zeta= e(a-c)\bar\zeta.\eeq For the above-given linear equation we look for a solution of the form $\zeta=pe^{\lambda t}$ where $p$ is a constant and the eigenvalues $\lambda$ are the solutions to the equation: \beqn\label{eig}\lambda-e(a-c)\int_{-r}^0K(s)e^{s\lambda} =0.\eeqn Here we consider the parameter values $b>d$ and $c>a$ resulting in snowdrift games which are classified as a type of social dilemma game (see, {\it e.g.,} \cite{aydogmus2020does}). For these parameter values, the replicator ODE has a unique stable interior equilibrium $e=\frac{b-d}{a-b-c+d}.$  Following \cite[Theorem 9]{ruan} we have the following result regarding the stability of this equilibrium under replicator DDE: \begin{proposition}
If \beq \int_0^r sK(s)\,ds<\frac{1}{e(c-a)}\eeq then the interior equilibrium $e$ of \eqref{delrep} is asymptotically stable.
\end{proposition}
 Here the quantity $\int_0^r sK(s)\,ds$ is generally called as the {\it average delay}. The above result implies that the interior Nash equilibrium is stable under the delayed replicator equations provided that the {\it average delay} is sufficiently small. 
 
 Note that \eqref{delrep} reduces the discrete delay replicator equations if Kernel $K(s)$ is the Dirac function $\delta(r-s).$ In particular, we obtain the following replicator equations: \beqn\label{ddelrep} \dot z(t)=z(t)\bigl(1- z(t)\bigr) \Bigl((a-b-c+d)z(t-r)+b-d \Bigr)\eeqn For this equation we have the following results regarding the stability of the interior equilibrium: \begin{proposition}\label{prop2}Let $R:=e^{-1}(c-a)^{-1}.$ Then we have
     \begin{itemize}
     \item[(i)] If $0\leq r< R\frac{\pi}{2}$ then $e$ is stable.
     \item[(ii)] If $r> R\frac{\pi}{2}$ then $e$ is unstable.
     \item[(iii)] If $r= R\frac{\pi}{2}$  a Hopf bifurcation occurs at $z = e$; that is, periodic solutions bifurcate from  $z = e.$ This implies periodic solutions exist for $r\geq R\frac{\pi}{2}$ and they are stable.
     \end{itemize}
\end{proposition} The proof of this result can be obtained by following the discussion by \citet[p. 481]{ruan}. 

This implies that the interior equilibrium of the replicator equations may loose its stability when the underlying game is a snowdrift game. On the other hand, both pure strategy equilibria $0$ and $1$ remain unstable. Then we have a unique interior equilibrium (which may or may not be stable under delayed replicator dynamics) and solutions are bounded away from the boundary of the $\Delta_2$ ($\ie$ 0 and 1). Hence both hypotheses of Theorem \ref{timeav} hold and the time average of $z$ converges to $e.$ This also implies by Corollary \ref{cor121} the existence of a constant time $\tau_\veps$ such that \beq
Pr\Biggl( \frac{1}{\tau N}\sum_{0\leq\sigma\leq\tau} \Zn({\sigma})\not\in B_\veps(\mathbf p)\Bigg)\leq 2d\exp\Bigl(-\veps^2 c_\tau N\Bigr)
\feq for any $\tau_\veps<\tau.$

\begin{figure}[ht]

 \centering
\subfigure[Simulations for $r=4$ and $N=1000$]{\label{chaos11}
\includegraphics[width=6.5cm, height=6.5cm]{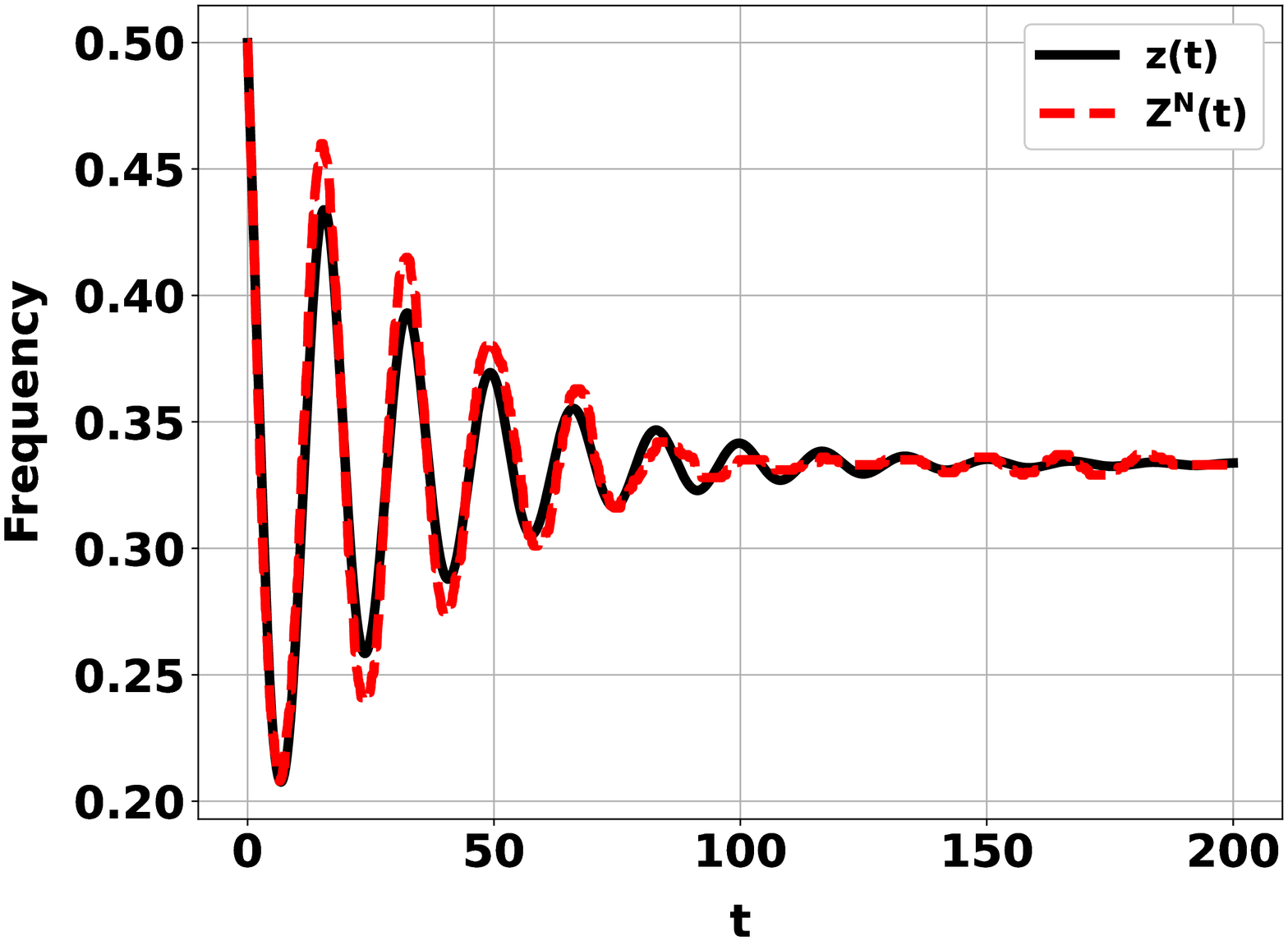}
}
\subfigure[Simulations for $r=4$ and $N=10000$]{\label{chaos21}
\includegraphics[width=6.5cm, height=6.5cm]{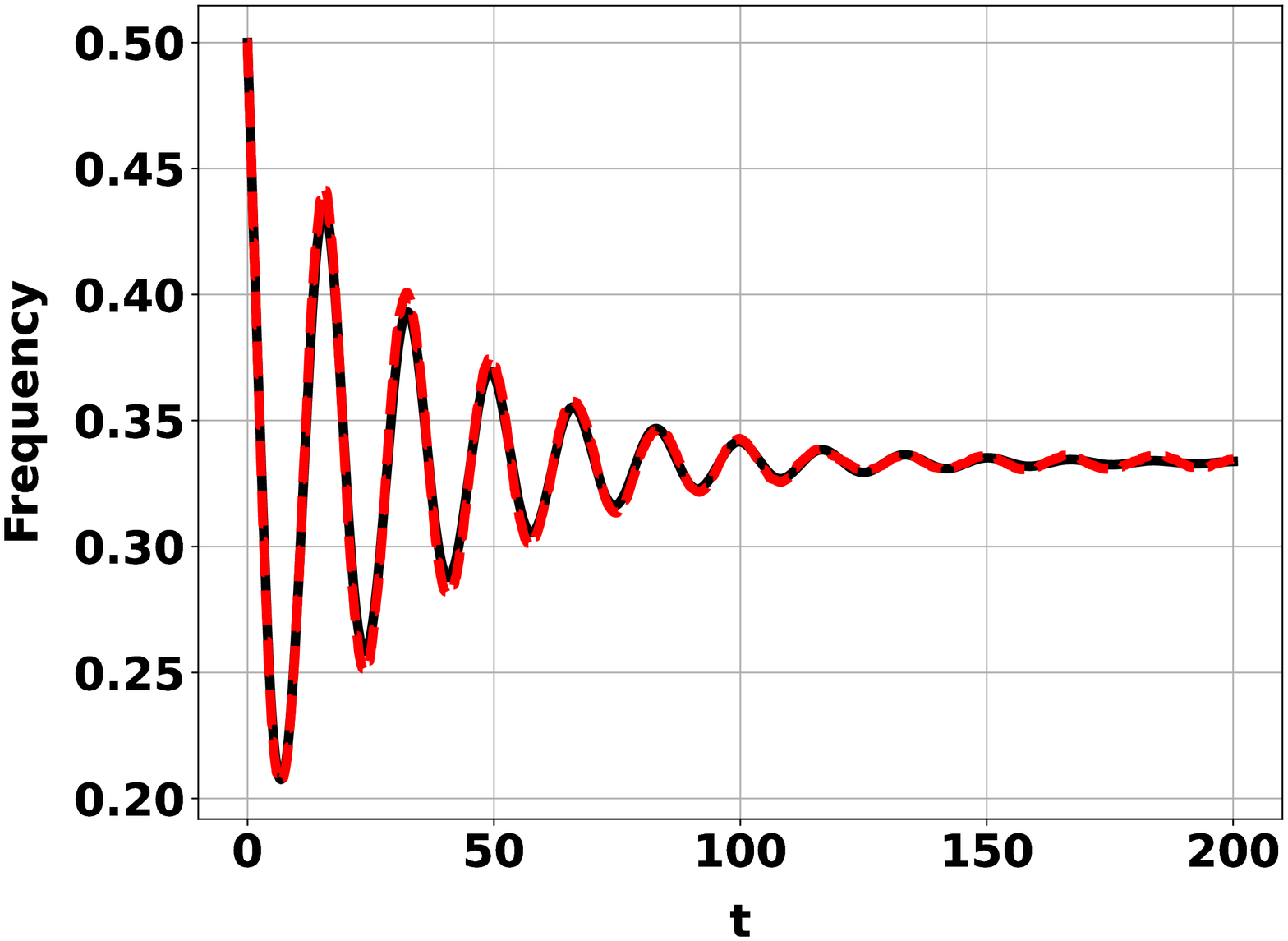}
}
\subfigure[Simulations for $r=5$ and $N=1000$]{\label{chaos31}
\includegraphics[width=6.5cm, height=6.5cm]{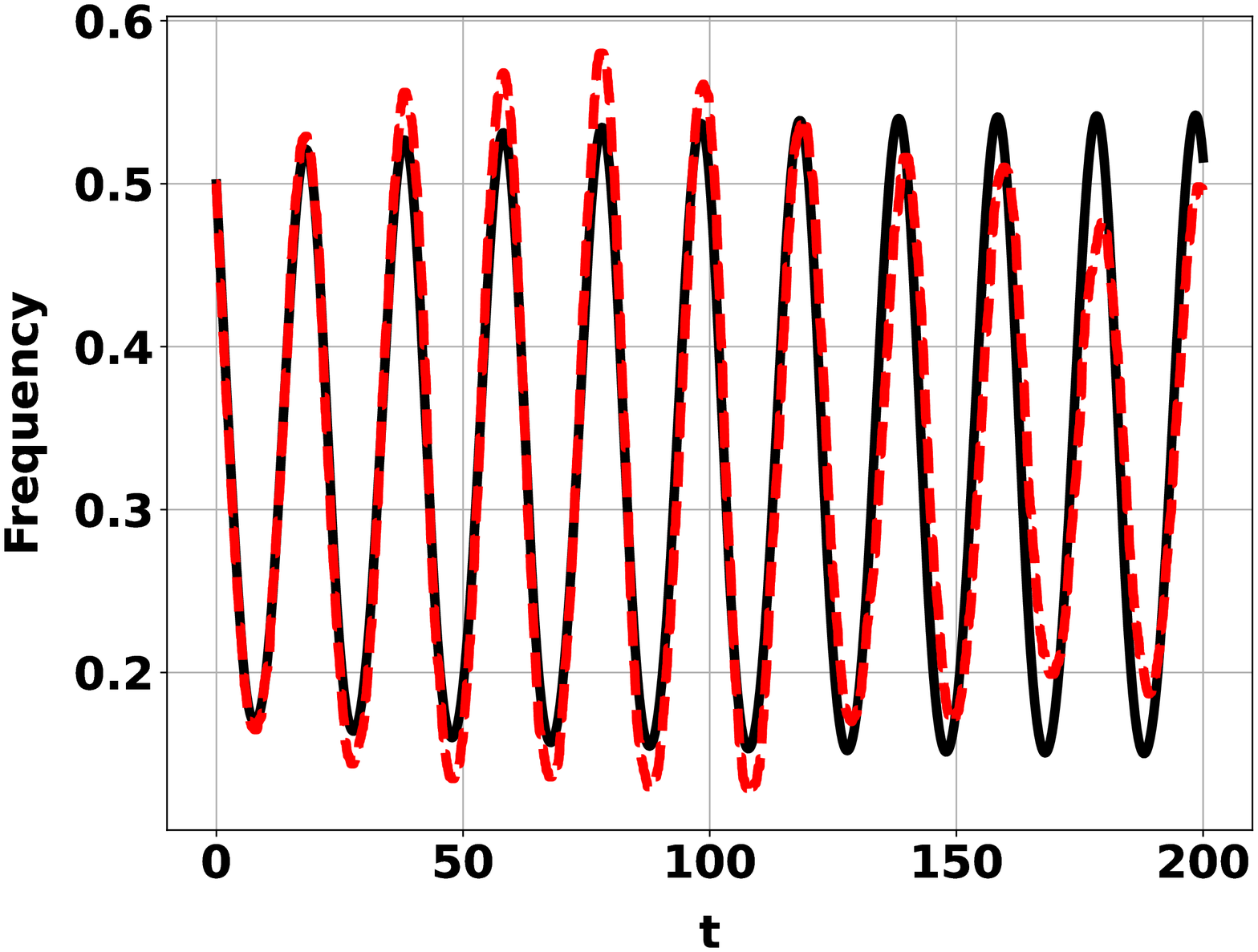}
}
\subfigure[Simulations for $r=5$ and $N=10000$]{\label{chaos311}
\includegraphics[width=6.5cm, height=6.5cm]{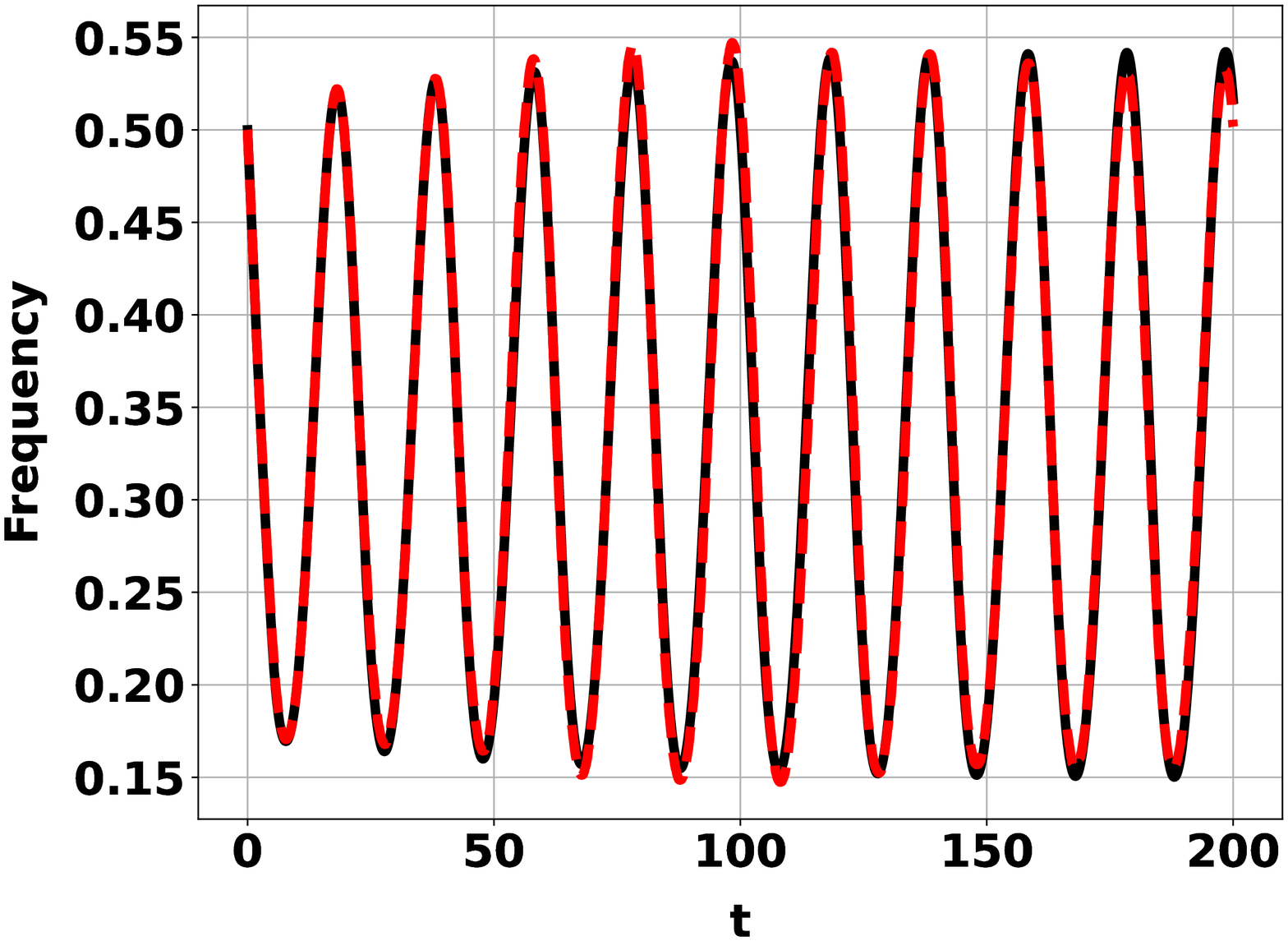}
}

\caption{  The cahnge of frequency of type A individuals for different values of the delay and different population sizes. }
\label{patternschaoticm}
\end{figure}

To verify above-given results numerically we considered a Hawk-Dove game with parameters $a=0.5,~b=0.5,~c=1.5$ and $d=0.$ Hence the unique interior Nash equilibrium of the game is given by $e=1/3.$ When the delayed replicator dynamics \eqref{ddelrep} is considered, $R=3$ (see Proposition \ref{prop2}). This implies that the critical value for the delay term $r$ is approximately $4.71.$ In our numerical simulations, we considered $r=4$ and $r=5.$ In Figures the upper figures( i.e. Figures \ref{chaos11}, \ref{chaos21}), we observe solution to the replicator equation approaches to the constant solution $e=1/ 3$ for $r=4$ which is less than the critical value. On the other hand, Figures \ref{chaos31} and \ref{chaos311} indicate that the periodic solutions emerge for $r=5$ that is larger than the critical value as noted in Proposition \ref{prop2}. In all of these figures, we compare numerically obtained trajectory of delay replicator equation \eqref{ddelrep} 
 with the trajectories obtained from non-Markovian process $Z$ for population sizes $N=1000$ (see Figures \ref{chaos11} and \ref{chaos31}) and $N=10000$ see Figures \ref{chaos21} and \ref{chaos311} up to time $200.$ As seen from the figures trajectories of deterministic and stochastic processes gets closer as the population size increases. The algorithms to simulate these deterministic and stochastic processes are given in Appendix \ref{algss}.
 \begin{figure}[ht]

 \centering

\includegraphics[width=6.5cm, height=6.5cm]{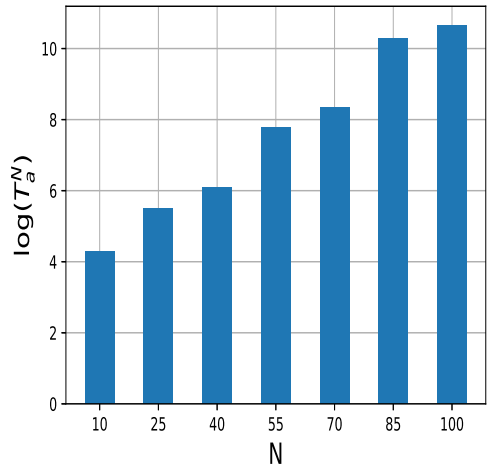}

\caption{Logarithm of Extinction times for different population sizes.}
\label{chaos41}
\end{figure}
 
 To verify the result regarding the exponentially increasing bound in population size for the fixation times given in Corollary \ref{corr}, we simulated the process for different population sizes ($N=10,25,40,55,70,85$ and $100$). The process is simulated until it reaches one of the absorbing states and this procedure is repeated 500 times for each population size. The average fixation times are found by taking the mean over these 500 samples and logarithm of fixation times $\log(T^N_a)$ versus the population size $N$ are illustrated in Figure \ref{chaos41}. The bar graph shows the existence of a linear relation between the logarithm of the fixation times and the population size. Hence this result verifies our theoretical finding given in Corollary \ref{corr}.

\section{Conclusion and Discussion}\label{conc}

We introduce a generalized class of birth-death processes that are
used to model the imitation dynamics in finite populations of interacting individuals by relaxing the assumption of memoryless agents. This class includes models in which individuals decide to change their strategies with a probability conditioned on the history of the process.  In particular, strategy update of an individual corresponds
to pairwise payoff comparison between individuals. Here we obtained {\it mean-field equations} of history dependent ({\it i.e.,} non-Markovian) processes and showed that these deterministic equations are good approximations to the finite population stochastic models in the sense that their trajectories stay arbitrarily close to each other up to a finite time $T$ with a probability approaching to one as the population size increases. Using this approximation result, we obtained two results regarding the fixation time of the process when the delayed replicator equations are bounded away from the boundary of the simplex: (1) The average fixation time increases exponentially with the population size and (2) the probability that the process never hits any of the absorbing states as the population size goes to infinity is one.  

 A well-known example of imitation dynamics is the replicator rule or replicator equations. In particular, it is shown that the {\it fluid limits} of so called replicator rules with the assumption of memoryless agents are replicator equations (see, {\it e.g.,} \cite{hofsigm,traulsen2005coevolutionary}). The replicator equations with discrete delays for specific two strategy and three strategy games are studied by \cite{alboszta2004stability,yi1997effect}, and \cite{wesson2016hopf}, respectively. Here, we showed that these equations are the fluid limits of the replicator rule taking the history of the process into account. In particular, we obtained replicator equations with discrete and distributed delays from the history dependent microscopic update rules for general $d$-strategy games and showed that the time averaging property of replicator ODEs (see, {\it e.g.,} \cite{hofsigm}) is also valid for the delayed replicator equations. 
 
Our model and analyses have limitations to keep the presentation of the paper simple. In particular, the results obtained here are not only valid for two player symmetric games. Multi-player asymmetric games can be considered following \cite{benaim}. Delayed replicator equations for three player games has been considered in \cite{bodnar2020three}. It is also possible to extend the model using infinite delays. Yet in a population of individuals imitating each other it is more realistic to assume that the agents use a recent and finite part of the history.   It is also possible to study the exit times from sets as done by \cite{benaim} instead of studying the fixation times. Yet such an extension requires to determine the basin of attractions of the deterministic equations which is in the space of continuous functions (denoted by $C$).  Since both determining such a subset of $C$ and attributing any meaning to exit time from this set  have no reasonable implications in the applied sciences, these results are not presented here. In addition, evolution of populations with strategy dependent time delays has been considered in \cite{mikekisz2021evolution}. Fluid limits of corresponding microscopic models of these processes may also be studied. The effect of delays on the dynamics of spatial models (see, e.g. \cite{aydogmus2017preservation,hwang2013deterministic,aydogmus2018discovering}) can also be studied. 

We extended the deterministic approximation results for Markov processes by \cite{benaim}. According to \cite{sandholm10}, these are the strongest approximation results in the literature. Our results can also be extended to continuous time non-Markovian processes following \cite{norris}. Here we considered a stochastic imitation (social) dynamics. Yet our results can be extended to study stochastic dynamics in population biology provided that a delay term is needed to model the growth of a population. The application areas include evolutionary games with biological-type time delay \citep{mikekisz2008evolutionary}, ecological models including Lotka-Volterrra equations \citep{kot2001elements,kuang} and epidemic models \citep{arino2006time}.

\section{Proofs}\label{prfs}
\subsection{Proof of Lemma \ref{exuniq}}\label{pexuniq}
 The existence and uniqueness of solutions to \eqref{delayde} on an interval $[-r,A]$ for some $A>0$ directly follows from Theorem 3.7 of \cite{smith2011}.

 Proposition 3.10 of \citet{smith2011} states that if a solution $\mathbf x:[-r,A)\to \rr^d$ to \eqref{delayde} with $0<A<\infty$ is noncontinuable then the solution must blow up as $t\to A.$ Hence, it is enough to show that the solutions to \eqref{delayde} are bounded which implies that solution to it can be extended to $[-r,\infty).$  When defining $p_{ij}$ above, we assumed that it satisfies $p_{ij}(\mathbf x_t)=0$ if $x_j(t)=0.$ By \eqref{veca}, this implies that $\dot x_i(t)=F_i(\mathbf x_t)\geq0$ for $x_i(t)=0$ and hence the non-negativity of the solutions.  Note also that $\sum_{i=1}^d F_i=0$ which implies the sum of all frequencies ($\sum_{i\in S}x_i$) remains constant all the time. Therefore, for any initial data $\boldsymbol\phi\in C$ satisfying $\boldsymbol\phi(\theta)\in \Delta_d,$ for all $\theta\in [-r,0]$ the solution $\mathbf x(t)$ stays in the simplex $\Delta_d$ and is bounded. This implies that a unique solution to the equation defined on $[-r,A)$ can be extended to a solution defined on $[-r,\infty).$\qed

\subsection{Proof of Theorem \ref{mainth}}\label{pmainth}
Here we follow the proof of Lemma 1 by \cite{benaim} and extend it to non-Markovian processes defined in Section \ref{stm}. We start by introducing the following notation  for the  {\it intorpolated process} \eqref{interpola}: \beq\mathbf{\Yn_t}(\theta)=\mathbf\Yn(t+\theta),~~ -r\leq \theta\leq 0.\eeq The following lemma is used in the proof. 
\begin{lemma}\label{lemma}
For aforementioned processes $\mathbf{\Yn_t}$ and $\overline{\mathbf{X}}^{\mathbf{(N)}}_{\mathbf t},$ we have \beq \big\|\mathbf{\Yn_t}-\overline{\mathbf{X}}^{\mathbf{(N)}}_{\mathbf t}\bigr\|_C\leq2\delta\eeq
\end{lemma}
\proof Recall from \eqref{notcs} that $\overline{\mathbf{X}}^{\mathbf{(N)}}(s):={\mathbf{X}}^{\mathbf{(N)}}(\sigma)$ for any $s\in [\sigma,\sigma+\delta).$ In addition, denote the extended set of transition times as follows: $\TT_e=\{-m\delta, (1-m)\delta,\cdots,-\delta,0,\delta,2\delta,\cdots\}.$ and define $I:=(t-r-\delta,t+\delta).$ Then we have\beq 
&&
\big\|\mathbf{\Yn_t}-\overline{\mathbf{X}}^{\mathbf{(N)}}_{\mathbf t}\bigr\|_C=\sup_{t-r\leq s\leq t}\bigl\|\mathbf \Yn(s)-\overline{\mathbf{X}}^{\mathbf{(N)}}(s)\bigr\|\\
&&
\qquad
\leq\max_{\sigma \in I\cap\TT_e} \Big\{\sup_{\sigma\leq s< \sigma+\delta}\bigl\|\mathbf \Yn(s)-{\mathbf{X}}^{\mathbf{(N)}}(s)\bigr\|~\Big\}\\
&&
\qquad
=\max_{\sigma \in I\cap\TT_e} \Big\{\sup_{\sigma\leq s< \sigma+\delta}\bigl\|\frac{s-\sigma}{\delta}\bigl(\mathbf\Xn (\sigma+\delta)-\mathbf\Xn (\sigma)\bigr)\bigr\|~\Big\}\\
&&
\qquad
\leq2\delta \eeq where one needs to employ Assumption $A$ if $t-r<0.$  This completes the proof.\qed

We denote the difference between the step taken by the stochastic process from $\tau$ to $\tau+\delta$ by\beq \mathbf U_{\tau}=N\bigl[\mathbf \Xn(\tau+\delta)-\mathbf\Xn(\tau)\bigr]-\mathbf F^m\bigl({\mathbf{X}}^{\mathbf{(N)}}_{\boldsymbol\tau}\bigr)\eeq where \beqn\label{discmf} F^m_i\bigl({\mathbf{X}}^{\mathbf{(N)}}_{\boldsymbol\tau}\bigr)=\sum_{k\neq i} p_{ik}^m\bigl({\mathbf{X}}^{\mathbf{(N)}}_{\boldsymbol\tau}\bigr) -\sum_{k\neq i}p_{ki}^m\bigl({\mathbf{X}}^{\mathbf{(N)}}_{\boldsymbol\tau}\bigr)\eeqn is the discrete version of \eqref{veca}.

Following \cite{benaim}, we give the following result which will be useful in proving Theorem \ref{mainth}.
\begin{lemma}\label{supm}
Let $\|\cdot\|_2$ denote the $L_2$ norm of a vector in $\rr^d.$ Then there exist a $\Gamma>0$ for which we have\beq E\Bigl[ e^{\langle \boldsymbol\Theta, \mathbf U_{\tau}\rangle}\big|\,\mathcal F_\tau\Big]\leq \exp\Big(\frac{\Gamma}{2}\|\boldsymbol\Theta\|_2^2\Big)\eeq for any $\boldsymbol\Theta\in\rr^d.$
\end{lemma}
\proof By \eqref{notds}, $F^m$ is a function of $m$ $d\mbox{-dimensional}$ vectors in $\Delta_d.$ Denote its maximum over these vectors as $\|F\|_2$ and take $\Gamma^2=(\sqrt{2}+\|F\|_2).$ The desired result follows from \cite[Lemma 3]{benaim}.\qed

In the following lines we give a proof of Theorem \ref{mainth}.

\proof Let $\mathbf U:\rr_+\to \rr^d$ be a map defined by 
$\mathbf U(t)=\mathbf U_\tau$ for $\tau\leq t<\tau+\delta.$  Suppose that $\tau$ be the largest element in $\TT$ satisfying $\tau\leq t$ then we clearly have the following equality: \beqn\label{equ1}
    \mathbf \Yn(t)-\boldsymbol{\phi}(0)=\frac{t-\tau}{\delta}\bigl[\mathbf\Xn (\tau&+&\delta)-\mathbf\Xn (\tau)\bigr]\nonumber\\
    &+& \sum_{\substack{ \sigma\in\TT\\ \delta\leq \sigma\leq t} } \bigl[\mathbf\Xn(\sigma)-\mathbf\Xn(\sigma-\delta)\bigr].
\eeqn where we used equality $\mathbf\Xn(0)=\boldsymbol\phi(0).$ For any $\sigma\in \TT$ we have \beqn\label{equ2} \mathbf\Xn(\sigma+\delta)&-&\mathbf\Xn(\sigma)=\int_{\sigma}^{\sigma+\delta} N\bigl[\mathbf\Xn(\sigma+\delta)-\mathbf\Xn(\sigma)\bigr] \,ds\nonumber\\
&=&\int_{\sigma}^{\sigma+\delta}\Big[ \mathbf F^m\big({\mathbf{X}}^{\mathbf{(N)}}_{\mathbf s}\big)+N\bigl[\mathbf\Xn(\sigma+\delta)-\mathbf\Xn(\sigma)\bigr]-\mathbf F^m\big({\mathbf{X}}^{\mathbf{(N)}}_{\mathbf s}\big)\Big]\,ds\nonumber\\
&=& \int_{\sigma}^{\sigma+\delta}\Big[ \mathbf F\big(\overline{\mathbf{X}}^{\mathbf{(N)}}_{\mathbf s}\big)+ \mathbf U(s)\Big]\,ds
\label{equ22}\eeqn where the last equality follows from \eqref{continuouspr} and \eqref{discmf}. Similarly it can be shown that \beqn\frac{t-\tau}{\delta}\bigl[\mathbf\Xn (\tau+\delta)-\mathbf\Xn (\tau)\bigr]=\int_{\tau}^{t}\Big[ \mathbf F\big(\overline{\mathbf{X}}^{\mathbf{(N)}}_{\mathbf s}\big)+ \mathbf U(s)\Big]\,ds.\label{equ3}\eeqn Using \eqref{equ22} and \eqref{equ3} in \eqref{equ1} gives the following equality:
\beqn\label{equmain}
    \mathbf \Yn(t)-\boldsymbol{\phi}(0)=\int_{0}^{t}\Big[ \mathbf F\big(\overline{\mathbf{X}}^{\mathbf{(N)}}_{\mathbf s}\big)+ \mathbf U(s)\Big]\,ds
\eeqn
In addition, we know that the continuous solution to \eqref{delayde} satisfies the following integral equation (see, $\eg,$ \cite{smith2011}): \beqn\label{dsint}
\mathbf x(t,\boldsymbol\phi)-\boldsymbol{\phi}(0)=\int_0^t \mathbf F\bigl(\mathbf x_s\big)\,ds.
\feqn
Therefore, by equations \eqref{equmain} and \eqref{dsint}, we have 
\beqn
    \mathbf \Yn(t)-\mathbf x(t,\boldsymbol\phi)&=&\int_{0}^{t}\Big[ \mathbf F\big(\overline{\mathbf{X}}^{\mathbf{(N)}}_{\mathbf s}\big) -\mathbf F\bigl(\mathbf x_s\big)+ \mathbf U(s)\Big]\,ds\nonumber\\
    &=&\int_{0}^{t}\Big[ \mathbf F\big(\overline{\mathbf{X}}^{\mathbf{(N)}}_{\mathbf s}\big)-\mathbf F\big(\mathbf{\Yn_s}\big)+\mathbf F\big(\mathbf{\Yn_s}\big) -\mathbf F\bigl(\mathbf x_s\big)+ \mathbf U(s)\Big]\,ds\nonumber
\eeqn Thus, for any $t\leq T,$ we have \beqn 
&& 
\| \mathbf \Yn(t)-\mathbf x(t,\boldsymbol\phi)\|=\nonumber
\\
&&
\qquad
=\Bigl\|\int_{0}^{t}\Big[ \mathbf F\big(\overline{\mathbf{X}}^{\mathbf{(N)}}_{\mathbf s}\big)-\mathbf F\big(\mathbf{\Yn_s}\big)+\mathbf F\big(\mathbf{\Yn_s}\big) -\mathbf F\bigl(\mathbf x_s\big)+ \mathbf U(s)\Big]\,ds\Bigr\|\nonumber\\
&&
\qquad
\leq 
\Bigl\|\int_0^t\mathbf U(s)\,ds\Bigr\| + K\int_{0}^{t}\Big[ \big\|\overline{\mathbf{X}}^{\mathbf{(N)}}_{\mathbf s}-\mathbf{\Yn_s}\big\|_C+\big\|\mathbf{\Yn_s} -\mathbf x_s\big\|_C\Big]\,ds\nonumber\\
&&
\qquad
\leq 
\Bigl\|\int_0^t\mathbf U(s)\,ds\Bigr\| + K\Bigg[2\delta T+\int_{0}^{t} \max_{-r\leq\mu\leq s}\big\|\mathbf{\Yn}(\mu) -\mathbf x (\mu,\boldsymbol\phi)\big\|\,ds\Bigg]\label{inek}
\eeqn where $K$ is the Lischitz constant for the map $F$ and the last inequality is due to Lemma \ref{lemma}. Now we let \beq v(s):= \max_{-r\leq\mu\leq s}\| \mathbf \Yn(\mu)-\mathbf x(\mu,\boldsymbol\phi)\|.\eeq Then inequality \eqref{inek} leads us to \beq v(t) \leq \psi (T) +K\Big[2\delta T+\int_0^t v(s)\,ds\Big]\eeq where \beq \psi (T)=\max_{t\in [0,T]}\Bigl\|\int_0^t\mathbf U(s)\,ds\Bigr\|.\eeq By Gr\"onwall's lemma, we have \beq D(T,\boldsymbol\phi)\leq v(T)\leq \bigl(\psi(T)+2K\delta T\bigr)e^{KT}\eeq
In particular, for $\delta\leq\frac{\veps}{4\,K\,T}e^{-K\,T},$ we have
\beqn
\nonumber 
 P\Bigl[D(T,\boldsymbol\phi) <\>\veps] 
\leq 
P\Bigl[\Psi(T)>\frac{\veps \exp(-K T)}{2}\Bigr]
\feqn where this equation is identical to equation (38) given of \cite{benaim}. In addition by Lemma \ref{supm}, we have the following super-martingale: \beq Z_\tau(\boldsymbol\Theta):=\exp\Big(\sum_{0\leq\sigma\leq\tau}\big\langle \boldsymbol\Theta,\mathbf U_{\sigma}\big\rangle-\frac{\Gamma}{2}\tau\delta\|\boldsymbol\Theta\|_2^2\Big).\eeq Following the proof \cite[Lemma 1]{benaim} we get \beq P\Bigl[\Psi(T)>\frac{\veps \exp(-K T)}{2}\Bigr]\leq2d\exp\Big(\veps^2\frac{e^{-2KT}}{8\delta\Gamma T}\Big).\eeq Hence, for $c=\displaystyle\frac{e^{-2KT}}{8 T(\sqrt{2}+\|\mathbf F\|_2)^{1/ 2}},$ we have the desired result.\qed

\subsection{Exit times from sets}\label{pth3}

The results concerning the absorption times depends on the exit time from sets. For a Borel subset $U\subset \Delta_d$ and an integer $N\in\nn,$ we denote 
\beq
T_U^N=\inf\Bigl\{\tau\geq 0:~\mathbf\Xn(\tau)\not\in U\Bigr\}.
\feq
Here  $T_U^N$ is the exit time of the stochastic process $\mathbf \Xn$ from the set $U.$ Regarding this quantity, we have the following result.
\begin{lemma}\label{lemmaexit}
 Let $\mathcal A\subset C$ be an attractor of the semiflow $\boldsymbol\Phi$ with the basin of attraction $B(\mathcal A)$ and suppose that $B\subset B(\mathcal A)$ is a compact set with $\mathbf Y^N_0\in B$ for all $N.$ Then there exist a constant $\alpha>0$ and a set $ U\subset \Delta_d$ such that\begin{itemize}
     \item[$\mathbf i$] $Pr\bigl[T_U^N\leq t\bigr]\leq2 (t+1) d\exp(-\alpha N)$
     \vspace{0.2cm}
     \item[$\mathbf{ii}$]$ E\bigl[T_U^N\bigr]\geq\frac{1}{4d}e^{\alpha N}-1$
      \end{itemize}
\end{lemma}

\proof ($\mathbf i$) Provided that $\mathcal A$ is an attractor of the semiflow, it is a compact set. Since $C$ is a subset of a normed vector space, it is path connected \citep[p. 200]{metric}. This implies that $\mathcal A$ is connected if and only if $B(\mathcal A)$ is connected \cite[p.32 ]{sell2013dynamics}. Since $B$ is also a compact set, we can find a bounded open neighborhood $V$ of $\mathcal A\cup B$ satisfying $\overline V\subset B(\mathcal A)$ where we used openness of $B(\mathcal A)$ (see \cite[Lemma 23.2]{sell2013dynamics}).  Since $\overline V$ is attracted to $\mathcal A,$ we have \beq S(t_0)\overline V\subset V\subset \overline V\subset B(\mathcal A)\eeq for some $t_0>\mathcal T.$ Hence for small enough $\veps,$ we have \beq N_\veps\bigl(S(t_0)\overline V\bigr)\subset V\subset N_\veps\bigl( V\bigr)\subset B(\mathcal A)\eeq
 Consider the time-$t_0$ map of the continuous time flow $S(t_0)\boldsymbol\phi$ for some $t_0>\mathcal T.$ For any $K\in \nn,$ we define the following  stochastic quantity
\beq  D_{K\,t_0}=\max_{0\leq k\leq K-1}  D\bigl(t_0, \mathbf Y^{\mathbf (N)}_{\mathbf{k\,t_0}}\bigr).\feq where $D(\cdot,\cdot)$ is as defined in \eqref{devia}. 
Then we have
\beqn
Pr[D_{K\,t_0}\geq\veps]&\leq& \sum_{k=0}^{K-1} Pr\biggl[D\bigl(t_0, \mathbf Y^{\mathbf (N)}_{\mathbf{k\,t_0}}\bigr) >\veps\biggr]\nonumber\\
&\leq& \sum_{k=0}^{K-1} E\biggl[\, Pr\bigl[ D\bigl(t_0, \mathbf Y^{\mathbf (N)}_{\mathbf{k\,t_0}}\bigr) >\veps|\mathbf Y^{\mathbf (N)}_{\mathbf{k\,t_0}} \bigl]\,\biggr]\nonumber\\
&\leq&\label{qqq1}2 K d\exp(-\veps^2 c_{t_0}N)
\feqn where the last inequality follows from Theorem \ref{mainth}.

Hence, $D_{K\,t_0}\leq\veps$ implies $\mathbf{\Yn_{t}}\in V$ for any $t\leq K\,t_0.$  Now we define the following set \beq U:=\bigcup \bigl\{\boldsymbol\phi(0)~|\,\phi\in V\biggr\}\eeq where $U\subset \Delta_d.$ Therefore, $\mathbf{\Yn_{t}}\in V$ for any $t\leq K\,t_0$ implies that $\mathbf{\Yn_}(t)\in U.$  
Then, for any $t\leq Kt_0,$ we have the following inequality \beq Pr\bigl[T_U^N \leq t\bigr] \leq Pr\bigl[D_{K\,t_0}\geq\veps\bigr]\leq2 K d\exp(-\veps^2 c_{t_0}N).
\feq Note that $ K$ can be chosen as $\lceil t/t_0\rceil.$ Note also that if $\mathcal T\leq1,$ $t_0$ can be chosen as $1.$ This implies that $ K\leq t+1.$  If, on the other hand, $\mathcal T>1$ then $ K\leq t+1.$ Thus, we have \beqn\label{ineq1}
Pr\bigl[T_U^N \leq t\bigr]\leq\label{qqq}2 (t+1) d\exp(-\veps^2 c_{t_0}N).
\feqn We get the desired result for $\alpha=\veps^2 c_{t_0}.$\\
($\mathbf{ii}$) The expected value can be computed via integrating the tail method as done by \cite[Lmma 4]{benaim}. Hence, we obtain $E[T_V^N]\geq \frac{1}{4d}\exp(\alpha N)-1.$\qed
\subsubsection{Proof of Corollary \ref{corr}}\label{pth31}
Since,  $\mathbf Y^N_0\in B(\mathcal A)$ we can find a compact set $B$ containing  $\mathbf Y^N_0$ and a bounded neighborhood $V$ of $B\cup\mathcal A.$  With these properties in hand, we know by Lemma \ref{lemmaexit} that $E\bigl[T_U^N\bigr]\geq\frac{1}{4d}e^{\alpha N}-1.$ If we show that $e_i\not\in U$ for all $i\in S$ we get $T_a^N\geq T_U^N.$ Suppose $\boldsymbol\phi(0)=\mathbf e_i$ for some $i\in S.$ Then $\boldsymbol\Phi(t,\boldsymbol \phi)$ is equal to  the constant function $\mathbf e_i\in C$ for all $t>0.$ This implies that $e_i\in\mathcal A$ which contradicts assumption $A.$ This completes the proof. \qed

\subsubsection{Proof of Corollary \ref{cor2}} \label{pcor2} After noting that $T_a^N\geq T_U^N$ for some set defined above, the proof basically follows from Borel-Cantelli Lemma and \eqref{ineq1} (for details see, \cite[lemma 2]{benaim}). 

\subsection{Delayed Replicator Equations}
\subsubsection{Proof of Theorem \ref{timeav}}\label{ptimeav}
 For any $i\in S,$ both equations \eqref{ddreplicator} and \eqref{distdreplicator} can be rewritten as follows: \beqn\label{newrep} \frac{d\log(x_i)}{dt}=e_i'A\mathbf y -\mathbf x'(t)A\mathbf y\eeqn where $\mathbf y$ is a vector valued function. Denote an accumulation point of time averages of $y_i$ by \beq z_i:=\lim_{T_m\to\infty}\frac{1}{T_m}\int_0^{T_m} y_i(t)\,dt\eeq for $i\in S.$ Then we say that $\mathbf z=(z_1,z_2,\cdots,z_d)\in \rr^d$ is a rest point of the replicator ODE (which can be obtained by taking $r=0$ in any of the equations \eqref{ddreplicator} or \eqref{distdreplicator}) if it satisfies the following conditions: \begin{itemize}
     \item[i-]  $\mathbf e_i'A\mathbf z=\mathbf e_j'A\mathbf z \mbox{  for any } i,j\in S$
     \item[ii-] $\mathbf z\in \Delta_d.$
     \end{itemize} To show (i) is true, we integrate both sides of \eqref{newrep} up to time $T$ to have \beq \log\bigl(x_i(T)\bigr)-\log\bigl(x_i(0)\bigr)=\sum_{j=1}^d a_{ij}\int_0^Ty_j(t)\,dt-\int_0^T\mathbf x'(t)A\mathbf y(t)\,dt\eeq for all $i\in S.$ By dividing both sides to $T$ and letting $T\to0,$ the left hand side of the equation vanishes. Then we get $\mathbf e_i A\mathbf z=\int_0^T\mathbf x'(t)A\mathbf y(t)\,dt$ for all $i\in S$ which completes the proof. 
     
     To show (ii) is correct we need to specify the function $\mathbf y.$\begin{itemize}
    \item Suppose that $\mathbf y(t)=\mathbf x(t-r)$ then \eqref{newrep} becomes \eqref{ddreplicator}. Moreover, we have \beq z_i:=\lim_{T_m\to\infty}\frac{1}{T_m}\int_0^{T_m} x_i(t-r)\,dt=\lim_{T_m\to\infty}\frac{1}{T_m}\int_0^{T_m} x_i(t)\,dt\eeq Hence $\mathbf z\in\Delta_d.$
     \item Now suppose that $\mathbf y(t)=\mathbf{\bar x}$ then \eqref{newrep} becomes \eqref{distdreplicator} and we have \beq z_i:&=&\lim_{T_m\to\infty}\frac{1}{T_m}\int_0^{T_m} \int_{-r}^0 x_i(t+s)\,d\mu(s)\,dt\\
     &=&\lim_{T_m\to\infty}\frac{1}{T_m}\int_0^{T_m} x_i(t)\,dt\eeq
    where we used Fubini's theorem and the dominated convergence theorem along with the fact that $\mu$ is a probability measure. This also implies that $\mathbf z\in\Delta_d.$
     \end{itemize} Above discussion leads us to the following facts:
     
     (a) If there is no interior equilibrium of the replicator ODE, any solution $x(t)$ of delayed equations \eqref{ddreplicator} or \eqref{distdreplicator} approaches to the boundary of the simplex $\partial\Delta_d$ asymptotically.
     
     (b) If there is a unique interior equilibrium $\mathbf p$ of the replicator ODE then any solution $x(t)$ with $x(0) \in \Delta_d\backslash\partial\Delta_d$ of delayed equations \eqref{ddreplicator} or \eqref{distdreplicator} approaches to $\mathbf p$ $\ie,$ $ p_i=z_i$ for all $i\in S.$
     
     \subsubsection{Proof of Corollary \ref{cor121}}\label{pc121}
First, consider the following definition\beq A_\tau\bigl(\mathbf x\bigr):=\frac{1}{\tau}\int_{0}^\tau \mathbf x(s)\,ds.\feq Second, note that \beq\sum_{0\leq\sigma\leq\tau} \mathbf{X}^{\mathbf{(N)}}({\sigma})=\sum_{1\leq\sigma\leq\tau}\int_{\sigma-\delta}^\sigma N\mathbf{\overline X}^{\mathbf{(N)}}(s)\,ds=N\int_0^\tau\mathbf{\overline X}^{\mathbf{(N)}}(s)\,ds\eeq Using this equality we have \beq A_\tau (\mathbf{\overline X^{(N)}})=\frac{1}{\tau N}\sum_{0\leq\sigma\leq\tau} \mathbf{X}^{\mathbf{(N)}}({\sigma})\eeq

Then, we have \beq \|A_\tau(\mathbf x)-A_\tau\bigl(\mathbf{\overline X^{(N)}}\bigr)\|&=&\biggl\|\frac{1}{\tau}\int_0^\tau \mathbf x(s)-\mathbf{\overline X^{(N)}}(s)\,ds\biggr\|\\
&\leq&\frac{1}{\tau}\int_0^\tau \max_{t\in[0,\tau]}\bigl\|\mathbf x(t)-\mathbf{\overline X^{(N)}}(t)\bigr\|\,ds\\
&\leq& \max_{t\in[0,\tau]}\bigl\|\mathbf x(t)-\mathbf{Y^{(N)}}(t)\bigr\|+\max_{t\in[0,\tau]}\bigl\|\mathbf \Yn(t)-\mathbf{\overline X^{(N)}}(t)\bigr\|
\feq By following the argument given in the proof of Lemma \ref{lemma}, we obtain the following inequality: \beq \|A_\tau(\mathbf x)-A_\tau\bigl(\mathbf{\overline X^{(N)}}\bigr)\|\leq D(\tau, N) +2\delta. \eeq Hence, we have \begin{eqnarray*} Pr\Bigl[ \|A_\tau(\mathbf x)-A_\tau\bigl(\mathbf{\overline X^{(N)}}\bigr\|>\epsilon\Bigr]\leq Pr[D(\tau,N)>\epsilon-2\delta ]\leq 2de^{-(\epsilon-2\delta)^2cN}\end{eqnarray*} by Theorem \ref{mainth}. By chosing $\epsilon=\frac{3\veps}{2}$ and assuming $\delta<\frac{\veps}{4}$ we have

\beq
 Pr\Bigl[ \|A_\tau(\mathbf x)-A_\tau\bigl(\mathbf{\overline X^{(N)}}\bigr)\|  >\frac{3\veps}{2} \Bigr]&\leq& 2d\exp\Bigl(-\veps^2 c_\tau N\Bigr).\feq Suppose that $ \|A_k(x_i)-A_k\bigl(X_i^N\bigr)\|_\infty  >\frac{3\veps}{2}.$ By Theorem \ref{timeav}, we can easily see that there exists a positive integer $\tau_\veps$ such that $A_\tau\bigl(x\bigr)\in B_{\veps/2}(\mathbf p)$ for any $\tau\geq \tau_\veps.$ These two facts imply that $\|A_k(X_i^N)-\mathbf p\|>\veps.$ Hence the desired result follows.

For the result regarding the absorption or fixation time $T_a,$ take $\veps=\min\{\|\mathbf p-\mathbf e_i\|/2~|\, i\in S\}.$ For such an $\veps>0,$ we have $\tau_\veps$ such that  \beq
Pr\Biggl( \frac{1}{\tau_0 N}\sum_{0\leq\sigma\leq\tau} \mathbf{X}^{\mathbf{(N)}}({\sigma})\in B_\veps(\mathbf p)\Bigg)\geq1- 2d\exp\Bigl(-\veps^2 c_{\tau_0} N\Bigr)
\feq for any fixed $\tau_0>\tau_\veps.$ For any $T\in \rr,$ define the variable $t_e:=\lceil T /\tau_0\rceil$ and observe that \beq
Pr\Biggl( \frac{1}{2T_a N}\sum_{0\leq\sigma\leq T_a} \mathbf{X}^{\mathbf{(N)}}({\sigma})\not\in B_\veps(\mathbf p)\Bigg)\leq 2d\exp\Bigl(-\veps^2 c_\tau N\Bigr)
\feq for any $\tau>\tau_\veps.$

\section*{Declarations}
\subsection*{Ethical Approval} Ethical approval is not required for this study.

\subsection*{Competing interests} Authors declare no competing interest.
\subsection*{Authors' contributions}
Both authors wrote and edited the paper.
\subsection*{Funding}
OA's research is supported by Fulbright Foundation and TUBITAK (The scientific and technological research council of Turkey) via the program 2219. OA is also thankful to
Arizona State University for its hospitality during a visit in which this work was carried out. YK's research is partially supported by NSF-DMS (Award Number 1716802\&2052820);  and The James S. McDonnell Foundation 21st Century Science Initiative in Studying Complex Systems Scholar Award (UHC Scholar Award 220020472). 

\subsection*{Availability of data and materials} Not applicable for this study.

\appendix
\section{Algorithms for simulations}\label{algss}
In the following we provide a basic finite difference algorithm to integrate the delay differential equation. For more stiff problems, Runge-Kutta methods can be employed or related libraries of any specific programming labguage may be used.\\
\begin{algorithm}
\caption{A finite difference algorithm to integrate the deterministic process}\label{alg:one}
\KwData{Predetermine parameters of payoff matrix $a,b,c,d,$ step size $\Delta t,$ delay parameter $r$ and a vector of initial data $y\in \mathbb R^{r/\Delta t}$}
\KwResult{A vector of the deterministic trajectory $z\in\mathbb R^{(200+r)/\Delta t}$}
$z[0:r\cdot N] \gets y$\;
\For{$i\gets r/\Delta t $ \KwTo $(200+r)/\Delta t $}{
    $z[i] \gets z[i-1]+\Delta t \bigg( z[i-1](1-z[i-1])\big((a-b-c+d)*z[i-r/\Delta t] +b-d\big)\bigg)$\; 
}
\end{algorithm}

The following algorithm, on the other hand, can be used to obtain a sample path of the stochastic process.
\begin{algorithm}
\caption{An algorithm for simulating the stochastic process}\label{alg:two}
\KwData{Predetermine parameters of payoff matrix $a,b,c,d,$ population size $N,$ delay parameter $r$ and a vector of initial data $Y\in \mathbb R^{r\cdot N}$}
\KwResult{A vector of the stochastic trajectory $Z\in\mathbb R^{(200+r)\cdot N}$}
$Z[0:r\cdot N] \gets Y$\;
\For{$i\gets r\cdot N$ \KwTo $(200+r)\cdot N$}{
    $fA \gets a Z[i-r\cdot N]+b(1-Z[i-r\cdot N])$\;
    $fB \gets c Z[i-r\cdot N]+d(1-Z[i-r\cdot N])$\;
    Generate a uniform random number $\rho \in [0,1]$\;
  \uIf{ $Z[i-1]\cdot(1-Z[i-1])\cdot(fA-fB)$ is greater than $\rho$}{
    $Z[i] \gets Z[i-1]+N^{-1}$\;
  }
  \uElseIf{$Z[i-1]\cdot(1-Z[i-1])\cdot(fB-fA)$ is greater than $\rho$}{
    $Z[i] \gets Z[i-1]-N^{-1}$\;
  }
  \Else{
    $Z[i] \gets Z[i-1]$\;
  }
}
\end{algorithm}

\bibliographystyle{plainnat}      
\bibliography{sample}   

\begin{thebibliography}{46}
\providecommand{\natexlab}[1]{#1}
\providecommand{\url}[1]{\texttt{#1}}
\expandafter\ifx\csname urlstyle\endcsname\relax
  \providecommand{\doi}[1]{doi: #1}\else
  \providecommand{\doi}{doi: \begingroup \urlstyle{rm}\Url}\fi

\bibitem[Alboszta and Miekisz(2004)]{alboszta2004stability}
Jan Alboszta and Jacek Miekisz.
\newblock Stability of evolutionarily stable strategies in discrete replicator
  dynamics with time delay.
\newblock \emph{Journal of theoretical biology}, 231\penalty0 (2):\penalty0
  175--179, 2004.

\bibitem[Arino and Van Den~Driessche(2006)]{arino2006time}
J~Arino and P~Van Den~Driessche.
\newblock Time delays in epidemic models.
\newblock In \emph{Delay differential equations and applications}, pages
  539--578. Springer, 2006.

\bibitem[Aydogmus(2016)]{aydogmus2016extinction}
Ozgur Aydogmus.
\newblock On extinction time of a generalized endemic chain-binomial model.
\newblock \emph{Mathematical biosciences}, 279:\penalty0 38--42, 2016.

\bibitem[Aydogmus(2018)]{aydogmus2018discovering}
Ozgur Aydogmus.
\newblock Discovering the effect of nonlocal payoff calculation on the stabilty
  of ess: Spatial patterns of hawk--dove game in metapopulations.
\newblock \emph{Journal of theoretical biology}, 442:\penalty0 87--97, 2018.

\bibitem[Aydogmus et~al.(2017)Aydogmus, Zhou, and
  Kang]{aydogmus2017preservation}
Ozgur Aydogmus, Wen Zhou, and Yun Kang.
\newblock On the preservation of cooperation in two-strategy games with
  nonlocal interactions.
\newblock \emph{Mathematical biosciences}, 285:\penalty0 25--42, 2017.

\bibitem[Aydogmus et~al.(2020)Aydogmus, Cagatay, and
  G{\"u}rpinar]{aydogmus2020does}
Ozgur Aydogmus, Hasan Cagatay, and Erkan G{\"u}rpinar.
\newblock Does social learning promote cooperation in social dilemmas?
\newblock \emph{Journal of Economic Interaction and Coordination}, 15\penalty0
  (3):\penalty0 633--648, 2020.

\bibitem[Bena{\"\i}m and Weibull(2003)]{benaim}
Michel Bena{\"\i}m and J{\"o}rgen~W Weibull.
\newblock Deterministic approximation of stochastic evolution in games.
\newblock \emph{Econometrica}, 71\penalty0 (3):\penalty0 873--903, 2003.

\bibitem[Binmore and Samuelson(1997)]{binmore97}
Ken Binmore and Larry Samuelson.
\newblock Muddling through: Noisy equilibrium selection.
\newblock \emph{journal of economic theory}, 74\penalty0 (2):\penalty0
  235--265, 1997.

\bibitem[Binmore et~al.(1995)Binmore, Samuelson, and Vaughan]{binmore}
Kenneth~G Binmore, Larry Samuelson, and Richard Vaughan.
\newblock Musical chairs: Modeling noisy evolution.
\newblock \emph{Games and economic behavior}, 11\penalty0 (1):\penalty0 1--35,
  1995.

\bibitem[Bodnar et~al.(2020)Bodnar, Miekisz, and Vardanyan]{bodnar2020three}
Marek Bodnar, Jacek Miekisz, and Raffi Vardanyan.
\newblock Three-player games with strategy-dependent time delays.
\newblock \emph{Dynamic Games and Applications}, 10\penalty0 (3):\penalty0
  664--675, 2020.

\bibitem[B{\"o}rgers and Sarin(1997)]{borgers}
Tilman B{\"o}rgers and Rajiv Sarin.
\newblock Learning through reinforcement and replicator dynamics.
\newblock \emph{Journal of economic theory}, 77\penalty0 (1):\penalty0 1--14,
  1997.

\bibitem[Boylan(1995)]{boylan}
Richard~T Boylan.
\newblock Continuous approximation of dynamical systems with randomly matched
  individuals.
\newblock \emph{Journal of Economic Theory}, 66\penalty0 (2):\penalty0
  615--625, 1995.

\bibitem[Broom and Krivan(2016)]{broom2018biology}
Mark Broom and Vlastimil Krivan.
\newblock Biology and evolutionary games.
\newblock In \emph{In: Basar, T., Zaccour, G. (eds) Handbook of Dynamic Game
  Theory}, pages 1--39. Springer, 2016.

\bibitem[Corradi and Sarin(2000)]{corradi}
Valentina Corradi and Rajiv Sarin.
\newblock Continuous approximations of stochastic evolutionary game dynamics.
\newblock \emph{Journal of Economic Theory}, 94\penalty0 (2):\penalty0
  163--191, 2000.

\bibitem[Cressman and Tao(2014)]{cressman2014replicator}
Ross Cressman and Yi~Tao.
\newblock The replicator equation and other game dynamics.
\newblock \emph{Proceedings of the National Academy of Sciences}, 111\penalty0
  (Supplement 3):\penalty0 10810--10817, 2014.

\bibitem[Darling and Norris(2008)]{norris}
RWR Darling and James~R Norris.
\newblock Differential equation approximations for markov chains.
\newblock \emph{Probab. Surv.}, 5:\penalty0 37--79, 2008.

\bibitem[Ethier and Kurtz(2009)]{ethier}
Stewart~N Ethier and Thomas~G Kurtz.
\newblock \emph{Markov processes: characterization and convergence}, volume
  282.
\newblock John Wiley \& Sons, 2009.

\bibitem[Ewens(2004)]{ewens2004mathematical}
Warren~John Ewens.
\newblock \emph{Mathematical population genetics: theoretical introduction},
  volume~1.
\newblock Springer, 2004.

\bibitem[Friedman(1998)]{friedman1998economic}
Daniel Friedman.
\newblock On economic applications of evolutionary game theory.
\newblock \emph{Journal of evolutionary economics}, 8\penalty0 (1):\penalty0
  15--43, 1998.

\bibitem[Fudenberg et~al.(2006)Fudenberg, Nowak, Taylor, and
  Imhof]{fudenberg2006evolutionary}
Drew Fudenberg, Martin~A Nowak, Christine Taylor, and Lorens~A Imhof.
\newblock Evolutionary game dynamics in finite populations with strong
  selection and weak mutation.
\newblock \emph{Theoretical population biology}, 70\penalty0 (3):\penalty0
  352--363, 2006.

\bibitem[Gintis et~al.(2003)Gintis, Bowles, Boyd, and
  Fehr]{gintis2003explaining}
Herbert Gintis, Samuel Bowles, Robert Boyd, and Ernst Fehr.
\newblock Explaining altruistic behavior in humans.
\newblock \emph{Evolution and human Behavior}, 24\penalty0 (3):\penalty0
  153--172, 2003.

\bibitem[Harms and Skyrms(2008)]{harms2008evolution}
William Harms and Brian Skyrms.
\newblock \emph{Evolution of moral norms}.
\newblock na, 2008.

\bibitem[Hofbauer and Sigmund(1998)]{hofsigm}
Josef Hofbauer and Karl Sigmund.
\newblock \emph{Evolutionary games and population dynamics}.
\newblock Cambridge university press, 1998.

\bibitem[Hwang et~al.(2013)Hwang, Katsoulakis, and
  Rey-Bellet]{hwang2013deterministic}
Sung-Ha Hwang, Markos Katsoulakis, and Luc Rey-Bellet.
\newblock Deterministic equations for stochastic spatial evolutionary games.
\newblock \emph{Theoretical Economics}, 8\penalty0 (3):\penalty0 829--874,
  2013.

\bibitem[Kot(2001)]{kot2001elements}
Mark Kot.
\newblock \emph{Elements of mathematical ecology}.
\newblock Cambridge University Press, 2001.

\bibitem[Kuang(1993)]{kuang}
Yang Kuang.
\newblock \emph{Delay differential equations with applications in population
  dynmaics}.
\newblock Academic Press, 1993.

\bibitem[Kurtz(1970)]{kurtz1970solutions}
Thomas~G Kurtz.
\newblock Solutions of ordinary differential equations as limits of pure jump
  markov processes.
\newblock \emph{Journal of applied Probability}, 7\penalty0 (1):\penalty0
  49--58, 1970.

\bibitem[Kurtz(1981)]{kurtz}
Thomas~G Kurtz.
\newblock \emph{Approximation of population processes}.
\newblock SIAM, 1981.

\bibitem[Mi{{e}}kisz(2008)]{mikekisz2008evolutionary}
Jacek Mi{{e}}kisz.
\newblock Evolutionary game theory and population dynamics.
\newblock In \emph{Multiscale problems in the life sciences}, pages 269--316.
  Springer, 2008.

\bibitem[Miekisz and Bodnar(2021)]{mikekisz2021evolution}
Jacek Miekisz and Marek Bodnar.
\newblock Evolution of populations with strategy-dependent time delays.
\newblock \emph{Physical Review E}, 103\penalty0 (1):\penalty0 012414, 2021.

\bibitem[Miekisz and Weso{\l}owski(2011)]{miekisz2011stochasticity}
Jacek Miekisz and Sergiusz Weso{\l}owski.
\newblock Stochasticity and time delays in evolutionary games.
\newblock \emph{Dynamic games and applications}, 1\penalty0 (3):\penalty0
  440--448, 2011.

\bibitem[Moreira et~al.(2012)Moreira, Pinheiro, Nunes, and
  Pacheco]{moreira2012evolutionary}
Jo{\~a}o~A Moreira, Flavio~L Pinheiro, Ana Nunes, and Jorge~M Pacheco.
\newblock Evolutionary dynamics of collective action when individual fitness
  derives from group decisions taken in the past.
\newblock \emph{Journal of theoretical biology}, 298:\penalty0 8--15, 2012.

\bibitem[O'Searcoid(2006)]{metric}
M{\'\i}che{\'a}l O'Searcoid.
\newblock \emph{Metric spaces}.
\newblock Springer Science \& Business Media, 2006.

\bibitem[Ruan(2006)]{ruan}
Shigui Ruan.
\newblock Delay differential equations in single species dynamics.
\newblock In \emph{Delay differential equations and applications}, pages
  477--517. Springer, 2006.

\bibitem[Sandholm(2010)]{sandholm10}
William~H Sandholm.
\newblock \emph{Population games and evolutionary dynamics}.
\newblock MIT press, 2010.

\bibitem[Sell and You(2013)]{sell2013dynamics}
George~R Sell and Yuncheng You.
\newblock \emph{Dynamics of evolutionary equations}, volume 143.
\newblock Springer Science \& Business Media, 2013.

\bibitem[Smith(2011)]{smith2011}
Hal~L Smith.
\newblock \emph{An introduction to delay differential equations with
  applications to the life sciences}, volume~57.
\newblock Springer New York, 2011.

\bibitem[Taylor et~al.(2004)Taylor, Fudenberg, Sasaki, and
  Nowak]{taylor2004evolutionary}
Christine Taylor, Drew Fudenberg, Akira Sasaki, and Martin~A Nowak.
\newblock Evolutionary game dynamics in finite populations.
\newblock \emph{Bulletin of mathematical biology}, 66\penalty0 (6):\penalty0
  1621--1644, 2004.

\bibitem[Traulsen and Hauert(2009)]{traulsen2009stochastic}
Arne Traulsen and Christoph Hauert.
\newblock Stochastic evolutionary game dynamics.
\newblock \emph{Reviews of nonlinear dynamics and complexity}, 2:\penalty0
  25--61, 2009.

\bibitem[Traulsen et~al.(2005)Traulsen, Claussen, and
  Hauert]{traulsen2005coevolutionary}
Arne Traulsen, Jens~Christian Claussen, and Christoph Hauert.
\newblock Coevolutionary dynamics: from finite to infinite populations.
\newblock \emph{Physical review letters}, 95\penalty0 (23):\penalty0 238701,
  2005.

\bibitem[Turner and Chao(1999)]{turner1999prisoner}
Paul~E Turner and Lin Chao.
\newblock Prisoner's dilemma in an rna virus.
\newblock \emph{Nature}, 398\penalty0 (6726):\penalty0 441--443, 1999.

\bibitem[Wakano and Aoki(2007)]{wakano2007social}
Joe~Yuichiro Wakano and Kenichi Aoki.
\newblock Do social learning and conformist bias coevolve? henrich and boyd
  revisited.
\newblock \emph{Theoretical population biology}, 72\penalty0 (4):\penalty0
  504--512, 2007.

\bibitem[Wang et~al.(2017)Wang, Yu, Kurokawa, and Tao]{wang2017imitation}
Shi-Chang Wang, Jie-Ru Yu, Shun Kurokawa, and Yi~Tao.
\newblock Imitation dynamics with time delay.
\newblock \emph{Journal of Theoretical Biology}, 420:\penalty0 8--11, 2017.

\bibitem[Weibull(1997)]{weibull1997evolutionary}
J{\"o}rgen~W Weibull.
\newblock \emph{Evolutionary game theory}.
\newblock MIT press, 1997.

\bibitem[Wesson and Rand(2016)]{wesson2016hopf}
Elizabeth Wesson and Richard Rand.
\newblock Hopf bifurcations in delayed rock--paper--scissors replicator
  dynamics.
\newblock \emph{Dynamic Games and Applications}, 6\penalty0 (1):\penalty0
  139--156, 2016.

\bibitem[Yi and Zuwang(1997)]{yi1997effect}
Tao Yi and Wang Zuwang.
\newblock Effect of time delay and evolutionarily stable strategy.
\newblock \emph{Journal of theoretical biology}, 187\penalty0 (1):\penalty0
  111--116, 1997.

\end{thebibliography}

%
%

\end{document}